\documentclass[11pt]{article}
\usepackage[pdftex, dvips]{graphicx} 
\usepackage{subfig}
\usepackage{epstopdf}
\usepackage{float}
\usepackage[english]{babel}
\usepackage{latexsym,amsmath,amssymb,amsfonts,cancel,dsfont}
\usepackage{bbm}
\usepackage{fancyhdr}
\usepackage[margin=2.5cm,top=3cm]{geometry}

\usepackage{cite}
\usepackage{csquotes}
\usepackage{xcolor}
\usepackage{caption}
\usepackage[colorlinks=true]{hyperref}

\newtheorem{defi}{Definition}[section]
\newtheorem{teo}[defi]{Theorem}
\newtheorem{prop}[defi]{Proposition}
\newtheorem{lem}[defi]{Lemma}

\newtheorem{remark}[defi]{Remark}

\newcommand{\dis}{\displaystyle}
\newcommand{\N}{\mathds{N}}
\newcommand{\Q}{\mathds{Q}}
\newcommand{\R}{\mathds{R}}

\newcommand{\jjnt}{\dis\int\!\!\!\!\int}
\newcommand{\la}{\lambda}

\newcommand{\va}{\varphi}

\newcommand{\black}{\color{black}}

\begin{document}

\title{ \huge 
\textbf{Some Inverse Problems for the Burgers Equation and Related Systems}}

\author{
J. Apraiz\thanks{Universidad del Pa\'is Vasco, Facultad de Ciencia y Tecnolog\'ia, Dpto.\ Matem\'aticas, Barrio Sarriena s/n 48940 Leioa (Bizkaia), Spain. E-mail: {\tt jone.apraiz@ehu.eus}.},
\ \ A. Doubova\thanks{Universidad de Sevilla, Dpto.\ EDAN e IMUS, Campus Reina Mercedes, 41012~Sevilla, Spain, E-mail: {\tt doubova@us.es}.},
\ \ E. Fern\'andez-Cara\thanks{Universidad de Sevilla, Dpto.\ EDAN e IMUS, Campus Reina Mercedes, 41012~Sevilla, Spain, E-mail: {\tt cara@us.es}.}, 
\ \ M. Yamamoto\thanks{University of Tokyo, Japan, E-mail:  myama@next.odn.ne.jp.}
}

\date{}

\maketitle

\begin{abstract}
   In this article we deal with one-dimensional inverse problems concerning the Burgers equation and some related nonlinear systems (involving heat effects and/or variable density). In these problems, the goal is to find the size of the spatial interval from some appropriate boundary observations of the solution. Depending on the properties of the initial and boundary data, we prove uniqueness and non-uniqueness results. In addition, we also solve some of these inverse problems numerically and compute approximations of the interval sizes.
\end{abstract}

\vspace*{0,5in}
 



\section{Introduction}\label{Sec-1}

   This paper deals with some inverse problems for nonlinear time-dependent PDEs in one spatial dimension.
    
   The analysis and solution of inverse problems of many kinds has recently increased a lot because of their relevance in many applications:
   elastography and medical imaging, seismology, potential theory, ion transport problems or chromatography, finances, etc.;
   see for instance~\cite{Borcea_et_al, Hanke, Richter}.
   The variety of inverse problems is huge in comparison with their direct analogs and many inverse problems coming from very classical and basic direct problems wait for theoretical and numerical research.
   Let us mention the monographs by~Bellassoued and~Yamamoto~\cite{BelYam3}, Isakov~\cite{Isakov}, Romanov~\cite{Romanov} and~Hasanov and~Romanov~\cite{HasRom}, where many theoretical and numerical aspects of inverse problems for partial differential equations are depicted.

   In this paper, we consider problems related to the identification of the size of the spatial interval where a time-dependent governing nonlinear equation must be satisfied.
   We will focus on the Burgers equation and some variants, satisfied for~$(x,t) \in (0,\ell) \times (0,T)$.
   We will assume that the equation is complemented with boundary and initial conditions corresponding to known data, respectively for $x \in \{0,\ell\}$ and $t = 0$.
   Then, we will try to determine the width $\ell$ of the spatial interval from some extra information, for instance given at $x = 0$.
   The main goals will be to establish or discard uniqueness and to compute approximations of the solutions to the inverse problems.
   Related questions have been analyzed recently for the linear heat and wave equations in~\cite{ACDFCY}.
       
   The plan is the following.
   In Section~\ref{Sec-2}, we consider the viscous Burgers equation under several different circumstances.
   Section~\ref{Sec-3} and~Section~\ref{Sec-4} respectively deal with the Burgers equation coupled to a heat equation and the variable density Burgers system.
   Finally, we present the results of some numerical experiments in~Section~\ref{Sec-5}.
   
\

   Throughout this paper, $\|\cdot\|$ and~$(\cdot\,,\cdot)$ will stand for the  usual $L^2$ norm and scalar product, respectively.
   In the particular case of the space $L^2(0,\ell)$, we will sometimes write $(\cdot\,,\cdot)_{\ell}$ in order to make explicit the length~$\ell$.
   The symbol $C$ will denote a generic positive constant.

\ 
\section{Some positive and negative results for the viscous Burgers equation}\label{Sec-2}

   Let us consider the following system for the Burgers equation:
   \begin{equation}\label{pb-Burgers}
\begin{cases}
u_t - u_{xx} + uu_x = 0, &0 < x < \ell,\ 0 < t < T, \\
u(0,t)=\eta(t),\ \ u(\ell,t)=0, & 0 < t < T, \\
u(x,0)=u_0(x), & 0<x<\ell.
\end{cases}
   \end{equation}
   
   The unknown $u = u(x,t)$ can be interpreted (for example) as the velocity of the particles of a homogeneous viscous fluid in a tube where the flow is allowed only lengthwise.
   It  can also be viewed as the car traffic density in a road in a simplified model, see for instance~\cite{Traffic}.
   
   The main inverse problem for~\eqref{pb-Burgers} is the following:

\
   
\noindent
\textbf{IP-1:} \textit{Fix $u_0 = u_0(x)$ and $\eta = \eta (t)$ in~\eqref{pb-Burgers} in appropriate spaces and assume that $\beta := u_x|_{x=0}$ is known.
   Then, find $\ell$.}

\

   We are first interested in proving uniqueness. More precisely, the following question is in order:
   
\

\noindent
\textbf{Uniqueness for IP-1:} \textit{Let $u^{\ell}$ and $u^L$ be the solutions to~\eqref{pb-Burgers} respectively associated to the spatial intervals $(0,\ell)$ and $(0,L)$.
   Assume that the corresponding observations $u_x^\ell(0,\cdot)$ and~$u_x^L(0,\cdot)$ coincide, that is,
   \begin{equation}\label{2a}
u^{\ell}_x(0,t) = u^{L}_x(0,t) \ \text{ in } \ (0,T).
   \end{equation}
   Then, do we have $\ell=L$?}

\

   In the sequel, we will provide some positive and negative answers to this question, depending on the kind of imposed boundary or initial data.

\bigskip
\subsection{The simplest cases: zero initial and/or boundary data}\label{SSec-2.1}

\subsubsection{ \textbf{Case I}: $\eta\not\equiv 0$ and $u_0\equiv 0$}\label{SSSec-2.1.1}

  If $u_0 \equiv 0$, we get uniqueness:
  
\begin{teo}\label{prop2.1-B}
   Assume that $0 < \ell \leq L$, $\eta \in L^\infty(0,T)$ satisfies $\eta\not\equiv 0$ and $u_0\equiv0$.
   Let $u^\ell$ and $u^L$ be the solutions to~\eqref{pb-Burgers} respectively corresponding to~$\ell$ and~$L$ and let us assume that, for some $M > 0$,
   \begin{equation}\label{2b}
|u_x^\ell(x,t)| \leq M \ \text{ in } \ (0,\ell) \times (0,T) \ \text{ and } \ |u_x^L(x,t)| \leq M \ \text{ in } \ (0,L) \times (0,T)
   \end{equation}
and~\eqref{2a} holds.
   Then, $\ell = L$.
\end{teo}

\noindent
{\bf Proof:}
   The proof is standard.
   It can be achieved by contradiction, assuming that $\ell< L$.
   Indeed, note that $u^\ell \in L^\infty((0,\ell) \times (0,T))$ and~$u^L \in L^\infty((0,L) \times (0,T))$.
   If we set $v := u^\ell - u^L$, one has
   \[
v_t - v_{xx} + v u_x^\ell + u^L v_x = 0 \ \text{ in } \ (0,\ell) \times (0,T)
   \]
and also $v(0,t) = 0$ and~$v_x(0,t) = 0$ in~$(0,T)$.
   Consequently, from the unique continuation property of the heat equation 
   (see~\cite{Saut-Scheurer}), we have $v = 0$ in~$(0,\ell) \times (0,T)$.
   This yields~$u^L(x,t) = 0$ in~$(\ell,L) \times (0,T)$ and then (again from unique continuation) $u^L \equiv 0$, which is an absurd. \hfill $\square$

\bigskip
\subsubsection{ \textbf{Case II}: $\eta\equiv 0$ and $u_0\not\equiv 0$}\label{SSSec-2.1.2} 

   Let us show that, as in the case of the linear heat equation (see~\cite{ACDFCY}), non-uniqueness holds in general.
   More precisely, a counter-example to uniqueness can be found.
   We will follow three steps:

\begin{itemize}

\item[1-] Using the Cole-Hopf transformation
   (named after J. D. Cole and E. Hopf's works~\cite{Cole} and~\cite{Hopf}, respectively), we will rewrite~\eqref{pb-Burgers} as a system for the heat equation.

\item[2-] Then, we will prove a result similar to~\cite[Proposition 2.1]{ACDFCY} and we will deduce non-uniqueness for the inverse problem corresponding to the heat equation with Neumann boundary conditions.
   
\item[3-] Finally, coming back to the original variables, we will be able to conclude.

\end{itemize}

   The Cole-Hopf transformation is given by 
   \begin{equation*}
\varphi(x,t)=e^{-\frac{1}{2}\int_0^xu(\xi,t)\,d\xi}
   \end{equation*}
or, equivalently,
   \begin{equation}\label{cole-hopf}
u(x,t)=-2\frac{\varphi_x(x,t)}{\varphi(x,t)} , \quad \varphi(0,t) \equiv 1.
   \end{equation}
   Using \eqref{cole-hopf}, the Burgers system~\eqref{pb-Burgers} can be rewritten in the form
   \begin{equation}\label{pb-heat}
\begin{cases}
\varphi_t - \varphi_{xx}= 0, &0 < x < \ell,\ 0 < t < T, \\
\varphi_x(0,t)=0,\ \ \varphi_x(\ell,t)=0, & 0 < t < T, \\
\varphi(x,0)=\varphi_0(x), & 0<x<\ell,
\end{cases}
   \end{equation}
where we have introduced $\varphi_0(x) := e^{-\frac{1}{2}\int_0^xu_0(\xi)\,d\xi}$.

   Let us denote by~$\lambda_n$ and~$\tilde\va_n$
   (resp.~$\mu_n$ and~$\tilde\psi_n$) the eigenvalues and eigenfunctions of the Neumann Laplacian in~$(0,\ell)$
   (resp.~$(0,L)$).
   Then,
   \begin{equation*}\label{eigenv}
\left\{ \begin{array}{l} \dis
\la_n := \frac{n^2\pi^2}{\ell^2}, \quad n \in \N\cup\{0\}, 
\\ \dis
\tilde{\va}_n(x) :=\begin{cases}\dis
\sqrt{\frac{2}{\ell}}\cos \left(\frac{n\pi x}{\ell}\right), \quad n \in \N,
\\ \dis
\frac{1}{\sqrt{\ell}},\quad n=0,\quad 0<x<\ell,
\end{cases}
\end{array}\right.
   \end{equation*}
and
   \begin{equation*}\label{eigenf}
\left\{ \begin{array}{l} \dis
\mu_n := \frac{n^2\pi^2}{L^2}, \quad n \in \N\cup\{0\},
\\ \dis
\tilde{\psi}_n(x) :=\begin{cases} \dis
\sqrt{\frac{2}{L}}\cos\left( \frac{n\pi x}{L}\right), \quad n \in \N,
\\ \dis
\frac{1}{\sqrt{L}},\quad n=0,\quad 0<x<L.
\end{cases}
\end{array}\right.
   \end{equation*}

   The solutions to~\eqref{pb-heat} corresponding to~$\ell$ and~$L$ can be defined for all~$t > 0$.
   They are respectively given by
   \begin{align}
\varphi^{\ell}(x,t)=& \sum_{n=0}^{\infty}  (\varphi_0, \tilde{\va}_n)_{\ell} \, e^{-\la_nt} \, \tilde{\va}_n(x) ,\quad 0<x<\ell,\, t>0\label{sol1}
\\ \noalign{\noindent\mbox{and}}
\varphi^{L}(x,t)=& \sum_{n=0}^{\infty} (\varphi_0, \tilde{\psi}_n)_{L} \, e^{-\mu_nt} \, \tilde{\psi}_n(x) ,\quad 0<x<L,\, t>0.\label{sol2}
   \end{align}

   Recall that these scalar products are respectively given by
   \begin{equation*}
(f,g)_{\ell} := \int^{\ell}_0 f(x)g(x) \,dx \quad\text{and}\quad
(f,g)_L := \int^L_0 f(x)g(x) \,dx.
   \end{equation*}

   For any set $K$, let us denote by~$\#K$ the cardinal of~$K$.
   Then, the following holds:
\bigskip

\begin{prop}\label{prop2.2}
   If $L/\ell \in \Q$, then there exist initial data $\varphi_0$ verifying 
   \begin{equation}\label{con1}
\#{ \{n : (\varphi_0, \tilde{\va}_n)_{\ell} \ne 0\} } = \#{ \{n : (\varphi_0, \tilde{\psi}_n)_{L} \ne 0\} } =1,
   \end{equation}
such that $\varphi_x^{\ell}(0,t) = \varphi_x^{L}(0,t)$ for all~$t>0$.
   Thus, we can have non-uniqueness with initial data $\varphi_0$ satisfying~\eqref{con1} 
even if $|L-\ell|$ is arbitrarily small.
\end{prop}

\noindent
\textbf{Proof}:
   Let $m_0, n_0 \in \N$ be given such that $n_0 < m_0$ and $\ell = n_0 L/m_0$, that is, $m_0/L = n_0/\ell$.
   Let us choose $k_1, n_1 \in \N$ such that $n_1 = k_1 m_0/n_0$.
   Note that
   \begin{equation*}
\la_{k_1} = \frac{k_1^2\pi^2}{\ell^2} = \frac{n_1^2\pi^2}{L^2}
= \mu_{n_1}
   \end{equation*}
and set
   \begin{equation}\label{fi0}
\varphi_0(x) := \cos \left(\frac{k_1\pi x}{\ell}\right) + a = \cos \left(\frac{n_1\pi x}{L}\right) + a, \quad x\in \R,
   \end{equation}
where $a$ is a real constant.

   The functions in~\eqref{sol1} and~\eqref{sol2} corresponding to this $\varphi_0$ are the following:
   \begin{equation}\label{solele}
\varphi^{\ell}(x,t) =a+ e^{-\frac{k_1^2\pi^2}{\ell^2}t}\cos{\left(\frac{k_1\pi}{\ell}x\right)}
   \end{equation}
and
   \begin{equation}\label{soleLe}
\varphi^{L}(x,t) =a+ e^{-\frac{n_1^2\pi^2}{L^2}t}\cos{\left(\frac{n_1\pi}{L}x\right)}.
   \end{equation}
   Consequently,
   \begin{equation*}
\varphi^{\ell}_x(0,t)=\varphi^{L}_x(0,t)=0.
   \end{equation*}
\hfill $\square$

\

   From~\eqref{cole-hopf}, \eqref{solele} and~\eqref{soleLe}, we get
   \begin{equation*}
u^{\ell}(x,t) = \frac{2k_1\pi}{\ell} \frac{e^{-\frac{k_1^2\pi^2}{\ell^2}t} \sin{\left(\frac{k_1\pi}{\ell}x\right)}}
{e^{-\frac{k_1^2\pi^2}{\ell^2}t}\cos{\left(\frac{k_1\pi}{\ell}x\right)} + a}
\ \ \text{ and } \ \
u^{L}(x,t) = \frac{2n_1\pi}{L} \frac{e^{-\frac{n_1^2\pi^2}{L^2}t} \sin{\left(\frac{n_1\pi}{L}x\right)}}
{e^{-\frac{n_1^2\pi^2}{L^2}t}\cos{\left(\frac{n_1\pi}{L}x\right)} + a} .
   \end{equation*}
   
   If $a$ is sufficiently large, these functions are well defined, solve the Burgers systems respectively in~$(0,\ell) \times (0,T)$ and~$(0,L) \times (0,T)$ for
   \begin{equation*}
u_0(x) =  \frac{2k_1\pi}{\ell} \frac{\sin{\left(\frac{k_1\pi}{\ell}x\right)}}{\cos{\left(\frac{k_1\pi}{\ell}x\right)} + a}
= \frac{2 n_1\pi}{L} \frac{\sin{\left(\frac{n_1\pi}{L}x\right)}}{\cos{\left(\frac{n_1\pi}{L}x\right)} + a}
   \end{equation*}
and, moreover, satisfy~\eqref{2b}.

   This ends the proof of non-uniqueness in this case. \hfill $\square$

\bigskip

\subsection{Results where $\eta (t)\not\equiv 0$ and $u_0(x)\not\equiv 0$}\label{SSec-2.2}

   In order to prove uniqueness when both $\eta$ and $u_0$ are nonzero
   (and~$\eta$ is sufficiently large), we need an auxiliary result that concerns traces of functions in~$H^2(0,\ell)$:

\begin{lem}\label{traceteo-2}
   Let $L_* > 0$ be given.
   Then
   \begin{equation*}
\left| \frac{df}{dx}(0) \right| \le \frac{C(L_*)}{\ell^{3/2}}
\,\| f\|_{H^2(0,\ell)}
   \end{equation*}
for any $f \in H^2(0,\ell)$ and any~$\ell$ with~$0 < \ell \leq L_*$.
\end{lem}

   The proof is elementary.
   It can be found in~\cite{ACDFCY}.
   
\begin{teo}\label{uniqeta-B}
   Assume that $0 < \ell \leq L \leq L_*$, $0 < T_0 < T$,
   \[
\begin{array}{c}
u_x^\ell(0,t) = u_x^L(0,t) \ \text{ in } \ (0,T), \ \ \| u_0 \|_{L^2(0,L)} \leq M_0, \\
\noalign{\smallskip}
|u_x^\ell(x,t)| \leq M \ \text{ in } \ (0,\ell) \times (T_0,T) \ \text{ and } \ |u_x^L(x,t)| \leq M \ \text{ in } \ (0,L) \times (T_0,T),
\end{array}
   \]
where $L_*$, $M_0$ and $M$ are some positive constants.
   There exists $\delta_0>0$ (only depending on~$L_*$, $T_0$, $T$, $M_0$ and~$M$) such that, if
   \begin{equation}\label{cond-eta}
\int^{T}_{T_0} | \eta(t) |^2 \,dt \geq \delta_0,
   \end{equation}
one necessarily has $\ell = L$.
\end{teo}

{\black
\noindent
\textbf{Proof:}
   In this proof, we will denote by~$A$ the one-dimensional Dirichlet Laplacian in~$(\ell,L)$, with associated eigenvalues $0 < \zeta_1 < \zeta_2 < \cdots$.
   
   Let us assume that~$\ell < L$.
   Then, arguing as in the proof of~Theorem~\ref{prop2.1-B}, we deduce that
   \begin{equation}\label{BC-uL}
u^{L}(\ell,t) = u^{L}(L,t) = 0 \ \text{ in } \ (0,T).
   \end{equation}
   Therefore, from well known energy estimates, one has
   \[
\| u^L(\cdot\,,t)\|_{L^2(\ell,L)} = \| e^{-tA} u^L(\cdot\,,0) \|_{L^2(\ell,L)} \leq M_0 e^{-\zeta_1 t}
\quad \forall t \in (T_0,T),
   \]
where $\zeta_1$ is the first eigenvalue of $A$, that is, $\zeta_1=\pi^2(L-\ell)^{-2}$.

\

   Let us put $u^L = u^h + z$ for $t \in (T_0,T)$, with $u^h(\cdot\,,t) := e^{-(t-T_0)A}u^L(\cdot\,,T_0)$.
   Then, we have:
   \[
\| u^h(\cdot\,,t)\|_{H^2(\ell,L)} \leq \frac{M_0}{T_0} e^{-\zeta_1 T_0} \ \text{ in } \ (T_0,T).
   \]
   On the other hand,
   \[
z(\cdot\,,t) = \int_{T_0}^t e^{-(t-s)A} (u^L u^L_x)(\cdot\,,s) \,ds
   \]
and the standard parabolic regularity estimates and the fact that $|u^L_x| \leq M$ yield:
   \[
\| z \|_{L^2(T_0,T;H^2(\ell,L))} \leq C(M) \| u^L \|_{L^2(T_0,T;L^2(\ell,L))} \leq C(T,M_0,M) e^{-\zeta_1 T_0}.
   \]
   Therefore,
   \[
\| u^L \|_{L^2(T_0,T;H^2(\ell,L))} \leq C(T,M_0,M)\left(1 + \frac{1}{T_0}\right) e^{-\zeta_1 T_0}
   \]
and, from~Lemma~\ref{traceteo-2}, we get:
   \begin{equation}\label{expi2-B}
\| u^{L}_x(\ell,\cdot) \|_{L^2(T_0,T)} 
\le \frac{C(L_*,T,M_0,M)}{(L-\ell)^{3/2}}\left(1 + \frac{1}{T_0}\right)
\exp\left( -\frac{\pi^2T_0}{(L-\ell)^2}\right).
   \end{equation}
   Maximizing the right hand side with respect to $L-\ell$, we obtain:
\begin{equation*}
\| u^{L}_x(\ell,\cdot) \|_{L^2(T_0,T)}  \leq \frac{1}{T_0^{3/4}}\left(1 + \frac{1}{T_0}\right) C(L_*,T,M_0,M) .
   \end{equation*}

   Now, we can continue exactly as in the proof of Theorem 2.7 in \cite{ACDFCY} and deduce that, if~$\delta_0$ is large enough, we get a contradiction. \hfill $\square$

\bigskip
\section{The Burgers equation with heat effects}\label{Sec-3}

   The system is now
   \begin{equation}\label{pbBH}
\begin{cases}
u_t - u_{xx} + uu_x = k \theta, &0<x<\ell,\, t>0,\\
\theta_t - \theta_{xx} + u \theta_x = 0, &0<x<\ell,\, t>0,\\
u(0,t)=\eta(t),\ \ u(\ell,t)=0, &t>0,\\
\theta(0,t) = \lambda(t), \ \ \theta(\ell,t) = 0, &t>0,\\
u(x,0)=u_0(x),\ \ \theta(x,0)=\theta_0(x),& 0<x<\ell.
\end{cases}
   \end{equation}
   Here, $k \in \R$ is given.
   
   As before, $u$ can be interpreted as the velocity of the fluid particles in a one-direction flow.
   This time, we assume that heat effects are important and, consequently, the evolution of a temperature $\theta = \theta(x,t)$ must also be taken into account.
   
   We will deal with the following inverse problem:
   
\

\noindent
\textbf{IP-2:} \textit{Fix $(u_0,\theta_0)$ and $(\eta,\lambda)$ in~\eqref{pbBH} in appropriate spaces and assume that $\beta := u_x|_{x=0}$ and~$\alpha := \theta_x|_{x=0}$ are known.
   Then, find $\ell$.}

   This is the uniqueness property we will analyze:

  \
  
\noindent
\textbf{Uniqueness for IP-2:} \textit{Let $(u^{\ell},\theta^{\ell})$ and $(u^{L},\theta^{L})$ be the solutions to~\eqref{pbBH} associated to the spatial intervals $(0,\ell)$ and $(0,L)$, respectively. Assume that the corresponding observations $(u_x^\ell(0,\cdot),\theta_x^\ell(0,\cdot))$ and~$(u_x^L(0,\cdot),\theta_x^L(0,\cdot))$ coincide, that is,
   \begin{equation}\label{14a}
u^{\ell}_x(0,t) = u^{L}_x(0,t) \ \text{ and } \ \theta^{\ell}_x(0,t) = \theta^{L}_x(0,t) \ \text{ in } \ (0,T).
   \end{equation}
   Then, do we have $\ell=L$?}
   
\

   If $(u_0,\theta_0)\equiv (0,0)$, we have again uniqueness:
  
\begin{teo}\label{prop-z}
   Assume that $0 < \ell \leq L < T$, $\eta, \lambda \in L^\infty(0,T)$ satisfy $(\eta,\lambda)\not\equiv (0,0)$ and $(u_0,\theta_0)\equiv(0,0)$.
   Let $(u^\ell,\theta^\ell)$ and $(u^L,\theta^L)$ be the solutions to~\eqref{pbBH} respectively corresponding to~$\ell$ and~$L$ and let us assume that~\eqref{2b} holds for some $M > 0$ and, furthermore, \eqref{14a} is satisfied.
   Then, $\ell = L$.
\end{teo}

   The proof is very similar to the proof of Theorem~\ref{prop2.1-B}.
   Thus, if we assume that $\ell < L$ and we set $v := u^\ell - u^L$ and~$\psi := \theta^\ell - \theta^L$, it is clear from unique continuation that $(v,\psi) = (0,0)$ in~$(0,\ell) \times (0,T)$.
   From energy estimates, we deduce that $(u^L,\theta^L) = (0,0)$ in~$(\ell,L) \times (0,T)$ and finally, again from unique continuation, $(u^L,\theta^L) \equiv (0,0)$, which is impossible. 
   
\

   On the other hand, it is obvious that any solution to~\eqref{pb-Burgers} is a particular solution to~\eqref{pbBH}, corresponding to~$\theta_0 \equiv 0$ and $\lambda \equiv 0$.
   Consequently, the counter-example considered in~Section~\ref{SSSec-2.1.2} is also a counter-example to uniqueness for~{\bf IP-2} when we allow $u_0$ to be nonzero.
   
   To our knowledge, it is unknown if a counter-example to uniqueness can also be found with $\theta_0 \not\equiv 0$.
   
\

   As before, we can deduce a uniqueness result for~\eqref{pbBH} for large $\eta$.
   More precisely, the following holds:

\begin{teo}\label{uniqeta-BS}
   Assume that $0 < \ell \leq L \leq L_*$, $0 < T_0 < T$,
   \[
\begin{array}{c}
\| (u_0,\theta_0) \|_{L^2(0,L)} \!\leq\! M_0, \ \ 
|u_x^\ell(x,t)| \!\leq\! M \ \text{ in } \ (0,\ell) \!\times\! (T_0,T), \ \ |u_x^L(x,t)| \!\leq\! M \ \text{ in } \ (0,L) \!\times\! (T_0,T)
\end{array}
   \]
and~\eqref{14a} holds.
   There exists $\delta_1>0$ (only depending on~$L_*$, $T_0$, $T$, $M_0$ and~$M$) such that, if
   \begin{equation}\label{cond-eta2}
\int^{T}_{T_0} | \eta(t) |^2 \,dt \geq \delta_1,
   \end{equation}
one necessarily has $\ell = L$.
\end{teo}

{\black
\noindent
\textbf{Proof:}
   It is similar to the proof of~Theorem~\ref{uniqeta-B}.
 
   Thus, let us assume that~$\ell < L$.
   As before, this implies 
   \[
u^{L}(\ell,t) = u^{L}(L,t) = 0 \ \text{ and } \ \theta^{L}(\ell,t) = \theta^{L}(L,t) = 0 \ \text{ in } \ (0,T).
   \]
   The following estimates for $(u^L,\theta^L)$ hold:
   \[
\| u^L(\cdot\,,T_0)\|_{L^2(\ell,L)} = M_0 e^{-\zeta_1 T_0}\ \text{and} \ \| \theta^L(\cdot\,,T_0)\|_{L^2(\ell,L)} = M_0 e^{-\zeta_1 T_0},
   \]
   \[
\| u^L \|_{L^2(T_0,T;L^2(\ell,L))} = C(T,M_0) e^{-\zeta_1 T_0} \ \text{ and the same hold for~$\theta^L$.}
   \]
   Let us put $u^L = w + z$, with $w(\cdot\,,t) := e^{-(t-T_0)A}u^L(\cdot\,,T_0)$.
   Then
   \[
\| w(\cdot\,,t)\|_{H^2(\ell,L)} \leq \frac{C}{T_0} \, e^{-\zeta_1 T_0} 
\ \text{ and } \ z(\cdot\,,t) = \int_{T_0}^t e^{ (t-s)A } \left(-u^L u_x^L + k \theta^L\right)(\cdot\,,s) \,ds
\ \text{ in } \ (T_0,T),
   \]
whence
   \[
\| z \|_{L^2(T_0,T;L^2(\ell,L))} \leq C \left[ M \| u^L \|_{L^2(T_0,T;L^2(\ell,L))} + k \| \theta^{L} \|_{L^2(T_0,T;L^2(\ell,L))} \right] .
   \]
   Consequently,
   \[
\| u^L \|_{L^2(T_0,T;H^2(\ell,L))} \leq C(T,L_*,M,M_0) \left(1+\frac{1}{T_0}\right) e^{-\zeta_1 T_0}
   \]
and
   \[
\| u_x^L(\ell,\cdot) \|_{L^2(T_0,T)} \leq \frac{C(T,L_*,M,M_0)}{(L-\ell)^{3/2}} \left(1+\frac{1}{T_0}\right) e^{-\frac{\pi^2 T_0}{(L-\ell)^2}}.
   \]

   At this point, we can continue as in the proof of~Theorem~\ref{uniqeta-B} and deduce that, for $\delta_1$ large enough, \eqref{cond-eta2} leads to a contradiction.
\hfill $\square$
}

\bigskip

   It is interesting to note that, in this result, the size of~$\lambda$
   (that is, $\theta|_{x=0}$) is not relevant at all.

\begin{remark}{\rm
   A simplified version of~\eqref{pbBH} can be obtained if we skip the transport terms.
   We find the linear system
   \begin{equation}\label{pb-1l}
\begin{cases}
u_t - u_{xx} = k \theta, &0<x<\ell,\, t>0,\\
\theta_t - \theta_{xx} = 0, &0<x<\ell,\, t>0,\\
u(0,t)=\eta(t),\ \ u(\ell,t)=0, &t>0,\\
\theta(0,t) = \lambda(t), \ \ \theta(\ell,t) = 0, &t>0,\\
u(x,0)=u_0(x),\ \ \theta(x,0)=\theta_0(x),& 0<x<\ell
\end{cases}
   \end{equation}
   It is not difficult to check that the assertions on uniqueness/nonuniqueness in Section~\ref{Sec-2} can be extended to this system with very similar (and in fact simpler) arguments.
\hfill $\square$
}
\end{remark}   

%
%
%
%
%
%
%
%

   Similar inverse problems can be considered for coupled Burgers-heat systems where a boundary temperature is observed.
   These are the following:
   \begin{equation}\label{pb-2}
\begin{cases}
u_t - u_{xx} + uu_x = k \theta, &0<x<\ell,\, t>0,\\
\theta_t - \theta_{xx} + u \theta_x = u_x^2, &0<x<\ell,\, t>0,\\
u(0,t) = \overline{u}(t), \ \ u(\ell,t) = 0, &t>0,\\
\theta_x(0,t) = \chi(t), \ \ \theta_x(\ell,t) = 0, &t>0,\\
u(x,0)=u_0(x),\ \ \theta(x,0)=\theta_0(x),& 0<x<\ell
\end{cases}
   \end{equation}
and
   \begin{equation}\label{pb-2l}
\begin{cases}
u_t - u_{xx} = k \theta, &0<x<\ell,\, t>0,\\
\theta_t - \theta_{xx} = 0, &0<x<\ell,\, t>0,\\
u(0,t) = u(\ell,t) = 0, &t>0,\\
\theta_x(0,t) = \chi(t), \ \ \theta_x(\ell,t) = 0, &t>0,\\
u(x,0)=u_0(x),\ \ \theta(x,0)=\theta_0(x),& 0<x<\ell.
\end{cases}
   \end{equation}
   
   \bigskip

   Now, the problems for~\eqref{pb-2} and~\eqref{pb-2l} are as follows: $(u_0,\theta_0)$, $\overline{u}$,  $\chi$ and the additional observations $\beta := u_x|_{x=0}$ and~$\zeta := \theta|_{x=0}$ are known and, again, we try to find $\ell$.

   The same questions above are in order.
   Results similar to Theorem~\ref{prop-z} and Theorem~\ref{uniqeta-BS} can be proved in this context.

\bigskip

\section{The case of the variable density Burgers equation}\label{Sec-4}

   This is more interesting, but also more difficult. We consider a non-homogeneous (or variable density) one-dimensional fluid, modeled as follows:
   \begin{equation}\label{pb-3}
\begin{cases}
\rho(u_t + uu_x) - u_{xx} = 0, &0<x<\ell,\, t>0,\\
\rho_t + u \rho_x = 0,  &0<x<\ell,\, t>0,\\
u(0,t)=\overline{u}(t),\ \ u(\ell,t)=0, &t>0, \\
\rho(0,t) = \overline{\rho}(t), &t \in \R_+ \cap \{ t : \overline{u}(t) > 0\}, \\
u(x,0)=u_0(x), \ \ \rho(x,0)=\rho_0(x), &0<x<\ell.
\end{cases}
   \end{equation}
 
   Of course, this can be viewed as a toy model for the variable density Navier-Stokes system.
   The corresponding inverse problem is the following:
   
\
   
\noindent
\textbf{IP-3:} \textit{Fix $(u_0,\rho_0)$ and~$(\overline{u},\overline{\rho})$ in~\eqref{pb-3} in appropriate spaces and assume that $\beta := u_x|_{x=0}$ and~$\gamma := \rho|_{x=0} 1_{\{t: \overline{u}(t) \leq 0 \}}$ are known.
   Then, find $\ell$.}
   
\

   This is the uniqueness question we are interested in:
  
\

\noindent
\textbf{Uniqueness for IP-3:} \textit{Let $(u^{\ell},\rho^{\ell})$ and $(u^{L},\rho^{L})$ be the solutions to~\eqref{pb-3} respectively associated to~$(0,\ell)$ and $(0,L)$.
   Assume that the corresponding $(u_x^\ell(0,\cdot),\rho^\ell(0,\cdot))$ and~$(u_x^L(0,\cdot),\rho^L(0,\cdot))$ coincide.
   Then, do we have $\ell=L$?}\\
}
\

\subsection{A result for zero initial data}

   When the initial data vanish, we have a positive uniqueness result for this problem:
   
\begin{teo}\label{prop-zz}
   Assume that $0 < \ell \leq L$, $T > 0$ and~$(u_0,\rho_0)$ and~$(\overline{u},\overline{\rho})$ satisfy
   \[
\left\{
\begin{array}{l} \dis
\overline{u}, \overline{\rho} \in L^\infty(0,T), \ \ \overline{u} \not\equiv 0, \ \ \overline{\rho} \geq 0,
\\ \noalign{\smallskip} \dis
u_0 \equiv 0, \ \ \rho_0 \in L^\infty(0,L), \ \ \rho_0 \geq a_0 > 0.
\end{array}
\right.
   \]
   Let $(u^\ell,\rho^\ell)$ and $(u^L,\rho^L)$ be the solutions to~\eqref{pb-3} for $0 < t < T$ respectively corresponding to~$\ell$ and~$L$. Let us assume that
$|u_t^\ell| + |u_x^\ell| + |\rho_x^\ell| \leq M$ and~$|u_t^L| + |u_x^L| + |\rho_x^L| \leq M$ respectively in~$(0,\ell) \times (0,T)$ and~$(0,L) \times (0,T)$ and $u^{\ell}_x(0,\cdot)=u^{L}_x(0,\cdot)$ and $\rho^{\ell}(0,\cdot)=\rho^{L}(0,\cdot)$. Then, $\ell = L$.
\end{teo}

   For the proof, we will use a unique continuation property satisfied by the solutions to systems of the form
   \begin{equation}\label{sys-aux}
\begin{cases}
a(x,t) v_t  - v_{xx} + b(x,t) v_x + c(x,t) v + d(x,t) p = 0,  & (x,t) \in Q,\\
p_t + m(x,t) p_x + r(x,t) v = 0, & (x,t) \in Q,\\
\end{cases}
   \end{equation}
where we assume that $Q := (0,\ell) \times (0,T)$,
   \begin{equation} \label{hyp-coeff}
a, b, c, d, m, r \in C^0(\overline{Q}) \ \text{ and } \ a \geq a_0 > 0 \ \text{ in } \ Q.
   \end{equation}
   More precisely, we have the following:

\begin{prop}\label{prop-Carleman}
   Assume that \eqref{hyp-coeff} is satisfied and~$(v,p)$ solves~\eqref{sys-aux}, with $v, v_x, v_{xx}, p, p_x \in C^0(\overline{Q})$.
   Also, assume that
  \begin{equation}\label{eq.vp}
\begin{cases}
v(0,t) = 0, \ \ v_x(0,t) = 0, \ \ p(0,t) = 0, &0<t<T, \\
v(x,0) = 0, \ \ p(x,0) = 0, &0<x<\ell.
\end{cases}
\end{equation}
   Then, one has $v \equiv 0$ and~$p \equiv 0$.
\end{prop}
   
   The proof of this Proposition relies on appropriate local Carleman estimates for the solutions to~\eqref{sys-aux} and is postponed to Section~\ref{SSec4.3}.
   
\

\noindent
{\bf Proof of Theorem~\ref{prop-zz}:}
   Note that $u^\ell \in L^\infty((0,\ell) \times (0,T))$ and~$u^L \in L^\infty((0,L) \times (0,T))$.
   If we set $v := u^\ell - u^L$ and~$p := \rho^\ell - \rho^L$, one has
   \[
\begin{cases}
\rho^\ell v_t  - v_{xx} + \rho^\ell v u^\ell_x + \rho^\ell u^L v_x + (u_t^L + u^L u_x^L) p = 0, 
&0<x<\ell,\, t>0,\\
p_t + u^L p_x + v \rho_x^\ell = 0,  &0<x<\ell,\, t>0,\\
v(0,t) = 0, \ \ v_x(0,t) = 0, \ \ p(0,t) = 0, &t>0, \\
v(x,0) = 0, \ \ p(x,0) = 0, &0<x<\ell.
\end{cases}
   \]
   
   Consequently, $v$ and~$p$ satisfies~\eqref{sys-aux} with $a = \rho^\ell$, $b = \rho^\ell u^L$, $c = \rho^\ell u_x^\ell$, $d = u_t^L + u^L u_x^L$, $m = u^L$ and $r = \rho_x^\ell$ and~\eqref{eq.vp}.
   
   In view of Proposition~\ref{prop-Carleman}, one has $v = 0$ and~$p = 0$ in~$(0,\ell) \times (0,T)$.
   This yields~$u^L(x,t) = 0$ in~$(\ell,L) \times (0,T)$.
   Since the equations satisfied by $u^L$ and $\rho^L$ also possess the unique continuation property, we find that $u^L \equiv 0$, which is impossible. \hfill $\square$

\

   It would be interesting to find nonzero initial data $(u_0,\rho_0)$ such that uniqueness fails.
   On the other hand, it would also be interesting to prove a result similar to~Theorem~\ref{uniqeta-BS} asserting that, if the boundary data are large enough
   (with respect to the other data in the system), uniqueness is satisfied.
   However, to our knowledge these questions are open.
   
\

   A still more complex situation is found when we deal with a variable density fluid where thermal effects are relevant.
   For example, we can consider the variable density Boussinesq-like system
   \begin{equation}\label{pb-4}
\begin{cases}
\rho(u_t + uu_x) - u_{xx} = \theta, &0<x<\ell,\, t>0,\\
\rho(\theta_t + u \theta_x ) - \theta_{xx} = u_x^2, &0<x<\ell,\, t>0,\\
\rho_t + u \rho_x = 0,  &0<x<\ell,\, t>0,\\
u(0,t)=\overline{u}(t),\ \ u(\ell,t)=0, &t>0, \\
\rho(0,t) = \overline{\rho}(t), &t \in \R_+ \cap \{ t : \overline{u}(t) > 0\}, \\
\theta_x(0,t) = \theta_x(\ell,t) = 0, &t>0,\\
\rho(x,0)=\rho_0(x),\ \ u(x,0)=u_0(x), \ \ \theta(x,0) = \theta_0(x), & 0<x<\ell.
\end{cases}
   \end{equation}

   This is the related inverse problem: $(u_0,\theta_0,\rho_0)$ and~$(\overline{u},\overline{\rho})$ are given and the additional observations $\beta := u_x|_{x=0}$ and $\zeta := \theta|_{x=0}$ are known for~$t \in (0,T)$ and we try to find $\ell$.

   A result similar to~Theorem~\ref{prop-zz} can also be established in this case.
   The details are left to the reader.
   
   \bigskip

\subsection{Proof of Proposition~\ref{prop-Carleman}}\label{SSec4.3}

   The proof of Proposition~\ref{prop-Carleman} can be obtained by combining two Carleman inequalities that can be deduced for the solutions to the first and the second equation in~\eqref{sys-aux}.
   The main steps are the following:
   
\begin{itemize}

\item To choose a suitable weight function
   (the same in both inequalities);
   
\item To argue as in~\cite{Ya1} and \cite{HIY} and deduce appropriate estimates for~$v$ and~$p$.

\item Finally, to add and eliminate all undesirable terms in the right hand side.

\end{itemize}

\medskip\noindent
\textbf{Step 1:} Let us first recall some known Carleman estimates for the solutions to equations like those in~\eqref{sys-aux}.

   Thus, assume that $a$, $b$ and~$c$ are as in~Proposition~\ref{prop-Carleman} and set $Lv := av_t - v_{xx} + b v_x + c v$ for any suitable~$v$.
   For any $\lambda>0$, $\beta>0$, $x_0 > \ell$, $\delta > 0$ and~$T > 0$
   (to be definitively fixed below), we take
   \begin{equation}\label{weight}
\varphi(x,t) := e^{\lambda \psi(x,t)}, \ \text{ with } \ \psi(x,t) := |x-x_0|^2 - \frac{2\delta\beta}{T} |t-T/2|.
   \end{equation}

   Note that $\varphi$ can be used in the proof of the Carleman inequality in~Theorem~2.1 in~\cite[Ch. 4]{Ya1}.
   Consequently, the following holds:

\bigskip
\begin{teo}\label{th.Carlemanp}
   There exists $\lambda_0>0$ with the following property:
   for any $\lambda \geq \lambda_0$, there exist constants $s_0=s_0(\lambda)>0$ and $C_0=C_0(\lambda)$ such that
   \begin{equation}\label{ineq.Carlemap}
\begin{array}{l}
\dis\jjnt_Q \left(\dfrac{1}{s\varphi} (|v_t|^2 + |v_{xx}|^2) + s\lambda^2 \varphi |v_x|^2 + s^3 \lambda^4 \varphi^3 |v|^2 \right) e^{2s\varphi}\, dx\,dt  
\\[4mm]
\phantom{--}\le C_0 \dis \left( \jjnt_Q |Lv|^2 e^{2s\varphi}\, dx\,dt  
+ \int_0^T 
\Big(s^3\lambda^3 \varphi^3 |v|^2 + s\lambda \varphi |v_x|^2 + |v_t|^2\Big) 
e^{2s\varphi} \, dt \Big |_{x=0,\ell} \right.
\\[4mm]
\phantom{--} +  \left. s^2 \lambda^2 e^{C_0\lambda}\dis \int_0^\ell \Big( |v|^2 + |v_x|^2 \Big ) e^{2s\varphi} \,dx\Big |_{t=0,T} \right)
\end{array}
   \end{equation}
for all $s\geq s_0$ and any $v\in H^{2,1}(Q)$. 
\end{teo}

\bigskip

   Now, let $m$ be as in~\eqref{hyp-coeff} and let us set $B := \varphi_t + m \varphi_x$ and $Ep := p_t + m p_x$ for any~$p$.
   We can also adapt the proof of the Carleman estimate for transport equations in~Proposition~2.1 in~\cite[Ch. 3]{Ya1} and deduce the following result: 

\bigskip
\begin{teo}\label{th.Carlemantr}
   Assume that $\min_{(x,t)\in\overline{Q}} |B(x,t)| \ge B_0>0$.
   Then, there exist constants $s_0 > 0$ and $C > 0$ such that 
   \begin{equation}\label{ineq.Carlematr}
\begin{array}{l}
\dis s^2 \jjnt_Q |p|^2 e^{2s\varphi}\, dx\,dt 
\leq C \dis  \jjnt_Q |Ep|^2 e^{2s\varphi}\, dx\,dt  
\\[4mm]
\phantom{\dis s^2 \jjnt_Q \Big |p|^2 }
+ s \dis\int_0^T  m B |p|^2 e^{2s\varphi} \, dt \Big |_{x=0}^{x=\ell}
+ s \dis\int_0^\ell B |p|^2 e^{2s\varphi}\, dx \Big |_{t=0}^{t=T} 
\end{array}
   \end{equation}
for all $s\geq s_0$ and any~$p\in H^{1}(Q)$. 
\end{teo}

\bigskip\noindent
\textbf{Step 2:} Let us assume that $t_0\in (0,T)$ and $\delta>0$ is such that $0<t_0-\delta < t_0 +\delta < T$ and let us set
   \[
Q_\delta :=(0,\ell) \times (t_0-\delta,t_0+\delta).
   \]  
   Let us introduce the new variable $\tilde{t}$ with $\tilde{t} =  t_0 -\delta + \dfrac{2\delta}{T} t$ and the new function $\tilde{\varphi}$ with 
   \[
\tilde{\varphi}(x,\tilde{t}) := e^{\lambda \tilde{\psi}(x,\tilde{t})} \ \text{ and } \ \tilde{\psi}(x,\tilde{t}) := \psi(x,\tilde{t}) \equiv |x-x_0|^2 - \beta |\tilde{t}-t_0|.
   \]
   Then, \eqref{ineq.Carlemap} can be rewritten as an estimate in~$Q_\delta$.
   By denoting $\tilde{t}$ (resp.~$\tilde{\varphi}$) again by~$t$ (resp.~$\varphi$), the following is found:

   \begin{equation}\label{ineq.CarlemapQd}
\begin{array}{l}
\dis\jjnt_{Q_\delta} \Big (\dfrac{1}{s\varphi} (|v_t|^2 + |v_{xx}|^2) + s\lambda^2 \varphi |v_x|^2 + s^3 \lambda^4 \varphi^3 |v|^2 \Big) e^{2s\varphi}\, dx\,dt  
\\[4mm]
\phantom{---} \leq C \left( \jjnt_Q |p|^2 e^{2s\varphi}\, dx\,dt + K_1 + K_2 \right),
\end{array}
   \end{equation}
where 

   \begin{equation}\label{K1}
\begin{array}{ll}
K_1  & := \dis  \int_{t_0-\delta}^{t_0+\delta} \Big(s^3\lambda^3 \varphi^3 |v|^2 + s\lambda \varphi |v_x|^2 + |v_t|^2\Big) e^{2s\varphi} \, dt \Big |_{x=0,\ell} 
\\[4mm]
& \leq
Cs^3\lambda^3 e^{C\lambda} \dis 
 \int_{t_0-\delta}^{t_0+\delta} \Big(|v(0,t)|^2 + |v_x(0,t)|^2 + |v_t(0,t)|^2\Big) e^{2s\varphi(0,t)} \, dt 
\\[4mm]
& \ \ \ +\ C s^3 \lambda^3 e^{C\lambda} M^2 \dis \int_{t_0-\delta}^{t_0+\delta} e^{2s\varphi(\ell,t)} \, dt
\end{array}
   \end{equation}
and
   \begin{equation}\label{K2}
\begin{array}{l}
K_2 := C s^2 \lambda^2 e^{C\lambda}\dis \int_0^\ell \Big( |v|^2 + |v_x|^2 \Big ) e^{2s\varphi}\, dx  \Big |_{t=t_0-\delta,t_0+\delta} \le 
Cs^2 \lambda^2 e^{C\lambda} M^2 e^{2s e^{\lambda (|x_0|^2 -\beta \delta)}}.
\end{array}
   \end{equation}

   On the other hand, the estimate (\ref{ineq.Carlematr}) applied to the second equation of~\eqref{sys-aux} in~$Q_\delta$ gives:
   
   \begin{equation*}
\begin{array}{l}
\dis s^2 \jjnt_{Q_\delta} |p|^2 e^{2s\varphi}\, dx\,dt  \leq C\dis  \jjnt_{Q_\delta} |v|^2 e^{2s\varphi}\, dx\,dt  
\\[4mm]
\phantom{\dis s^2 \jjnt_{Q_\delta} \Big |p|^2}
+ s \dis\int_{t_0-\delta}^{t_0+\delta}  m B |p|^2 e^{2s\varphi} \Big |_{x=0}^{x=\ell} \, dt
+ s \dis\int_0^\ell B |p|^2 e^{2s\varphi}\, dx  \Big|_{t = t_0-\delta}^{t = t_0+\delta}
\end{array}
   \end{equation*}
and we find that
   \begin{equation}\label{ineq.CarlematrQd}
\begin{array}{l}
\dis s^2 \jjnt_{Q_\delta} \Big |p|^2 e^{2s\varphi}\, dx\,dt \leq C\dis \jjnt_{Q_\delta} |v|^2 e^{2s\varphi}\, dx\,dt  + R_1 + R_2,
\end{array}
   \end{equation}
where 
   \begin{equation}\label{R1}
\begin{array}{c}
R_1  : =  C s e^{C\lambda} M^2 \dis\int_{t_0-\delta}^{t_0+\delta}  |p|^2 e^{2s\varphi}\, dt \Big |_{x=0}^{x=\ell} 
\\[4mm]
\le Cs e^{C\lambda} M^2 \dis \int_{t_0-\delta}^{t_0+\delta}|p(0,t)|^2 e^{2s\varphi(0,t)}\, dt 
+ Cs e^{C\lambda} M^4 \dis \int_{t_0-\delta}^{t_0+\delta} e^{2s\varphi(\ell,t)}\, dt 
\end{array}
   \end{equation}
and
   \begin{equation}\label{R2}
\begin{array}{l}
R_2 : = Cs e^{C\lambda} M \dis\int_0^\ell |p|^2 e^{2s\varphi}\, dx  \Big|_{t = t_0-\delta}^{t = t_0+\delta} \le C s e^{C\lambda} M ^3 e^{2s e^{\lambda(|x_0|^2 -\beta \delta)}} .
\end{array}
   \end{equation}

   In~\eqref{K1}, \eqref{K2}, \eqref{R1} and~\eqref{R2}, we have used that $|v|+|v_x|+|v_t| + |p| \leq M$ in $\overline{Q}$.
   It is not restrictive to assume that $M \geq 1$.

\bigskip\noindent
\textbf{Step 3:} After adding~\eqref{ineq.CarlemapQd} and~\eqref{ineq.CarlematrQd}, if we take into account the estimates of the $K_i$ and $R_i$ and the data and observations, assuming that $s$ and~$\lambda$ are sufficiently large, we find:

   \begin{equation}\label{Carleman_add}
\begin{array}{l}
\dis\jjnt_{Q_\delta} \Big (\dfrac{1}{s\varphi} (|v_t|^2 + |v_{xx}|^2) + s\lambda^2 \varphi |v_x|^2 + s^3 \lambda^4 \varphi^3 |v|^2 \Big) e^{2s\varphi}\, dx\,dt  
+ \dis s^2 \jjnt_{Q_\delta} \Big |p|^2 e^{2s\varphi}\, dx\,dt  
\\[4mm]
\phantom{---}\leq  C s^3 \lambda^3 e^{C\lambda} M^2 \dis\int_{t_0-\delta}^{t_0+\delta} \Big(|v(0,t)|^2 + |v_x(0,t)|^2 + |v_t(0,t)|^2 + |p(0,t)|^2\Big) e^{2s\varphi(0,t)} \, dt 
 \\[4mm]
\phantom{----} + C s^3 \lambda^3 e^{C\lambda} M^4 \dis\int_{t_0-\delta}^{t_0+\delta} e^{2s\varphi(\ell,t)} \, dt
+ C s^2 \lambda^2 e^{C\lambda}  M^3 e^{2s e^{\lambda(|x_0|^2 - \beta \delta)}}\\
\phantom{---} = C s^3 \lambda^3 e^{C\lambda} M^4 \dis\int_{t_0-\delta}^{t_0+\delta} e^{2s\varphi(\ell,t)} \, dt
+ C s^2 \lambda^2 e^{C\lambda}  M^3 e^{2s e^{\lambda(|x_0|^2 - \beta \delta)}}.
\end{array}
   \end{equation}

\bigskip

   Now, we argue as follows:
   
\begin{itemize}

\item First, we fix $\lambda > 0$ such that \eqref{Carleman_add} holds and choose $x_0$, $t_0$ and~$\delta$ as before and~$\varepsilon \in (0,\ell)$.

\item Then, we take $\beta > 0$ large enough and such that $\beta \delta /2 > \ell x_0 + \ell^2$.

\item Finally, we choose $\kappa \in (0,\delta/2)$ such that $\beta\kappa < 2 \varepsilon (x_0-\ell) + \varepsilon^2$.

\end{itemize}

   With these constants $\varepsilon$ and~$\kappa$, one has
   
   \begin{equation}\label{good-ineq}
|x \!-\! x_0|^2 \!-\! \beta|t \!-\! t_0| \!\geq\! \mu \!:=\! |x_0 \!-\! \ell \!+\! \varepsilon|^2 \!-\! \beta\kappa \!>\! 
\max(|x_0 \!-\! \ell|^2, |x_0|^2 \!-\! \beta\delta)
   \end{equation}

\noindent  
for all $(x,t) \in (0,\ell - \varepsilon) \times (t_0 - \kappa,t_0 + \kappa)$.
   Taking into account~\eqref{eq.vp}, we deduce from~\eqref{Carleman_add} that

   \begin{equation}\label{Carleman_casi}
\begin{array}{l}
\dis\jjnt_{(0,\ell - \varepsilon) \times (t_0 - \kappa,t_0 + \kappa)} \Big (s \lambda^4 |v|^2 + |p|^2 \Big) \, dx\,dt  
\\[4mm]
\phantom{---}\leq  2\delta C s \lambda^3 e^{C\lambda} M^4 e^{2s (e^{\lambda |x_0-\ell|^2} - e^{\lambda \mu})}
+ C s \lambda^2 e^{C\lambda}  M^3 e^{2s (e^{\lambda(|x_0|^2 - \beta \delta)} - e^{\lambda \mu})}
\\[4mm]
\phantom{---}\leq C_* s \left( e^{2s (e^{\lambda |x_0-\ell|^2} - e^{\lambda \mu})} + e^{2s (e^{\lambda(|x_0|^2 - \beta \delta)} - e^{\lambda \mu})} \right),
\end{array}
   \end{equation}
   
   \noindent
where $C_*$ depends on~$M$, $\delta$ and~$\lambda$ but is independent of~$s$.
   But, in view of~\eqref{good-ineq}, this right hand side goes to zero as $s \to +\infty$.
   Consequently, $v(x,t) = 0$ and~$p(x,t) = 0$ in~$(0,\ell - \varepsilon) \times (t_0 - \kappa,t_0 + \kappa)$.
   
   Since $\varepsilon$ and~$\kappa$ are arbitrarily small and~$t_0$ is arbitrary in~$(0,T)$, $v \equiv 0$ and~$p \equiv 0$ and the proof is achieved. \hfill $\square$

\bigskip
\section{Some numerical results}\label{Sec-5}

   In this section, we will perform  some numerical experiments for the previous inverse problems.
   We will carry out the reconstruction of the unknown length through the resolution of some appropriate extremal problems.
   This strategy has been applied in some previous papers of the authors for other similar problems, see~\cite{DouFe}, \cite{DouFe2} and~\cite{DouFe3}.
   The results of the numerical tests that follow will serve to illustrate the theoretical results in the previous sections. 

\subsection{Inverse problems for the Burgers equation}

   We deal with the following

\bigskip
\noindent
\textbf{Reformulation of IP-1:} \textit{Given $T > 0$, $\eta = \eta(t)$, $u_0 = u_0(x)$ and~$\beta = \beta(t)$, find $\ell \in (\ell_0,\ell_1)$ such that 
   \begin{equation}\label{pb.opt}
J_1(\ell) \leq J_1(\ell') \quad \forall\, \ell'\in (\ell_0,\ell_1) ,
   \end{equation}
where $J$ is given by
   \begin{equation}\label{eq.J}
J_1(\ell) := \dfrac{1}{2}\displaystyle\int_0^T |\beta(t) - u^\ell_x (0,t)|^2\, dt.
   \end{equation}
   Here, $u^\ell$ is the {\it state,} i.e.\ the solution to~\eqref{pb-Burgers} corresponding to the length~$\ell$.}

\bigskip
   
   Three different situations will be analyzed for the Burgers equation.
   In the first two cases, we will check that uniqueness holds: zero initial data and nonzero initial data and sufficiently large~$\eta$.
   In the third case we will consider a non-uniqueness situation corresponding to some nonzero initial data and ``small'' $\eta$ and we will study the behavior of the numerical  algorithm.
   To this purpose
   (and also in the experiences in the following sections), we will implement the~\texttt{fmincon} function from the MatLab Optimization Toolbox using the~\texttt{active-set} minimization algorithm. 

\

\noindent
\textbf{Case 1.1: Burgers equation with $u_0=0$ and $\eta\neq 0$.}

   We take $T=5$,  $\eta(t) = 5\sin^3 t$ in~$(0,T)$ and~$u_0(x) \equiv 0$.
   Starting from $L_{i} = 3$, our goal is to recover the desired value of the length $L_d= 2$.

\begin{figure}[h!]
\begin{minipage}[t]{0.49\linewidth}
\noindent
\includegraphics[width=\linewidth]{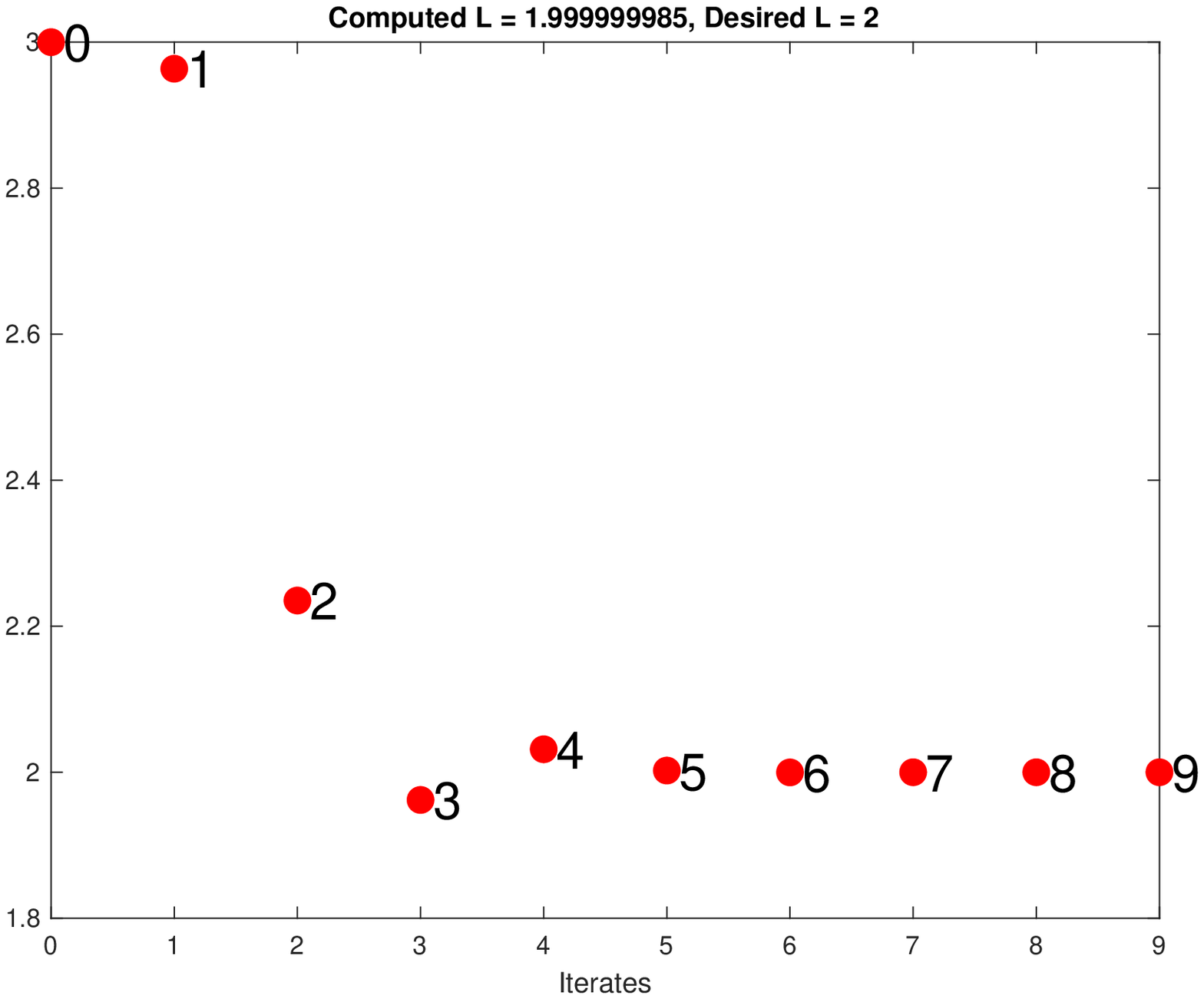}
\caption{Burgers equation, $u_0=0$ and $\eta\neq 0$. The iterates in \texttt{active-set} algorithm.}
\label{Burgers_iter_case1}
\end{minipage}
\hfill
\begin{minipage}[t]{0.49\linewidth}
\noindent
\includegraphics[width=\linewidth]{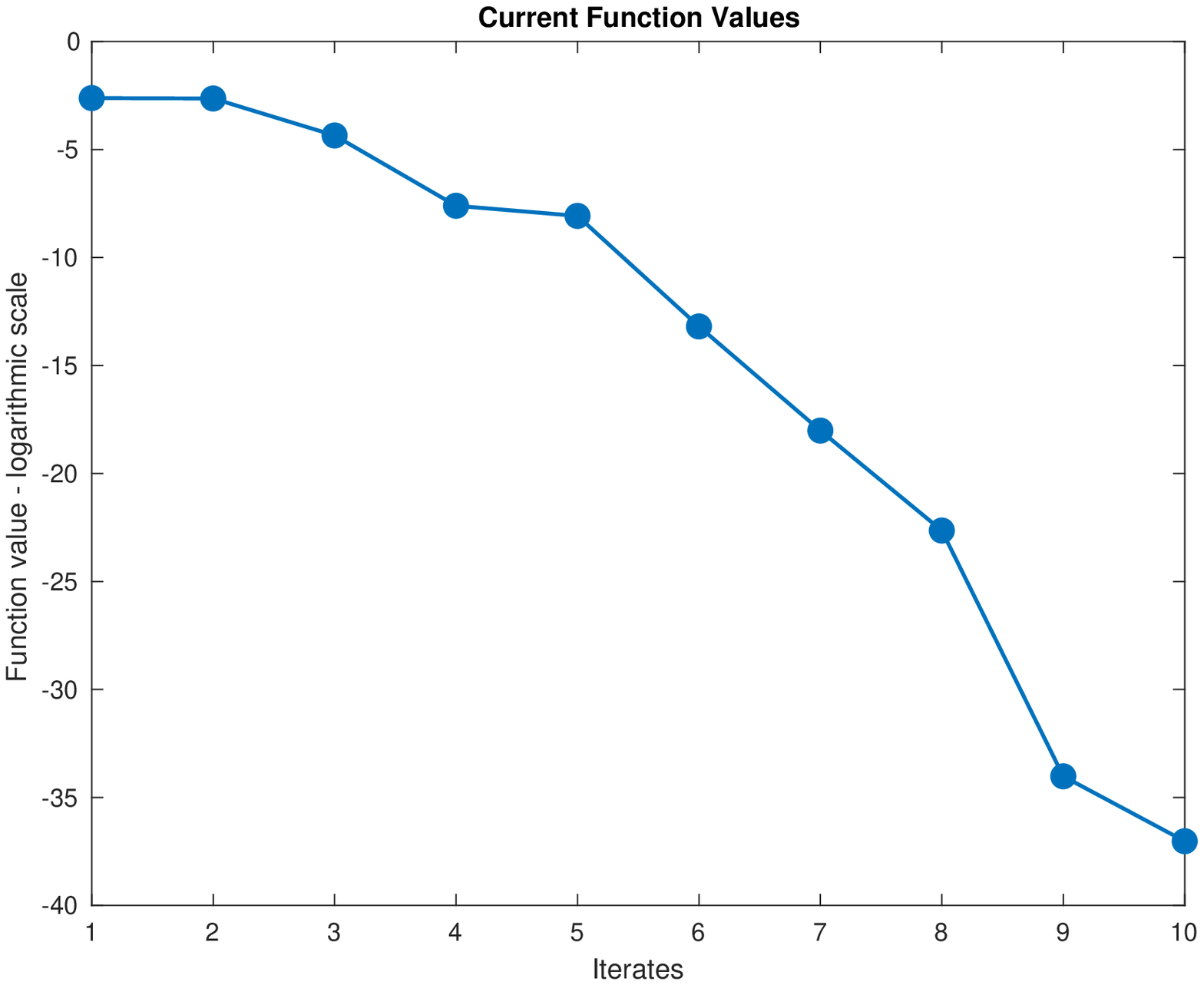}
\caption{Burgers equation, $u_0=0$ and $\eta\neq 0$. Evolution of the cost.}
\label{Burgers_J_case1}
\end{minipage}
\end{figure}

   The results of this numerical experiments can be seen in~Table~\ref{tab:burgers_case1}, where the effect of random noise on the target are shown.
   The computed length is denoted by~$L_c$.
   The corresponding solution to~\eqref{pb.opt}--\eqref{eq.J} is displayed in~Figure~\ref{Burgers_Sol_case1}.
   The evolution of the iterates and the cost in the minimization process in the absence of random noise appear in~Figures~\ref{Burgers_iter_case1} and~\ref{Burgers_J_case1}, respectively.
   

%
\begin{table}[h!]
\begin{minipage}[t]{0.46\linewidth}
\centering
\renewcommand{\arraystretch}{1.2}
\medskip
\caption{Burgers equation, $u_0=0$ and $\eta\neq 0$. Results with random noise in the target
(desired length: $L_d = 2$). \\
{ }
{ }}
\begin{tabular}{|cccc|} \hline
\%  noise  & Cost  & Iterates &  Computed $L_c$ \\
\hline
1\%           &   1.e-3         &       12       &   1.997140631 \\ \hline
0.1\%        &    1.e-5        &      15        &   1.999169558\\ \hline
0.01\%      &   1.e-7       &       11         &   1.999912907\\\hline
0.001\%    &   1.e-9       &     10         &   2.000021375\\ \hline
0\%           &   1.e-17       &     9         &    1.999999985  \\  \hline
\end{tabular}
\label{tab:burgers_case1}
\end{minipage}
\hfill
\begin{minipage}[t]{0.54\linewidth}
\vspace{0cm}
\includegraphics[width=\linewidth]{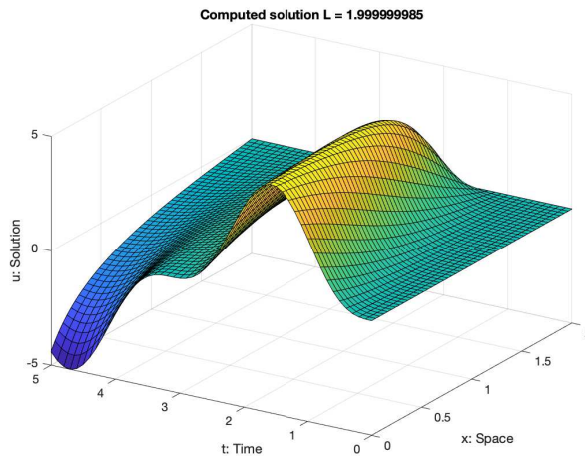}
\captionof{figure}{Burgers equation with $u_0=0$ and $\eta\neq 0$. The computed solution.}
\label{Burgers_Sol_case1}
\end{minipage}
\end{table}
%


\noindent 
\textbf{Case 1.2: Burgers equation with $u_0\neq 0$ and large $\eta$.}

   We take $T=5$,  $\eta(t) = 5 \sin(t)^3$ in~$(0,T)$ and~$u_0(x) \equiv 3 x (2-x)$.
   Now, starting from $L_{i} = 2.4$, the target value that we want to recover is $L_d= 2$.
   
\begin{figure}[h!]
\begin{minipage}[t]{0.49\linewidth}
\noindent
\includegraphics[width=\linewidth]{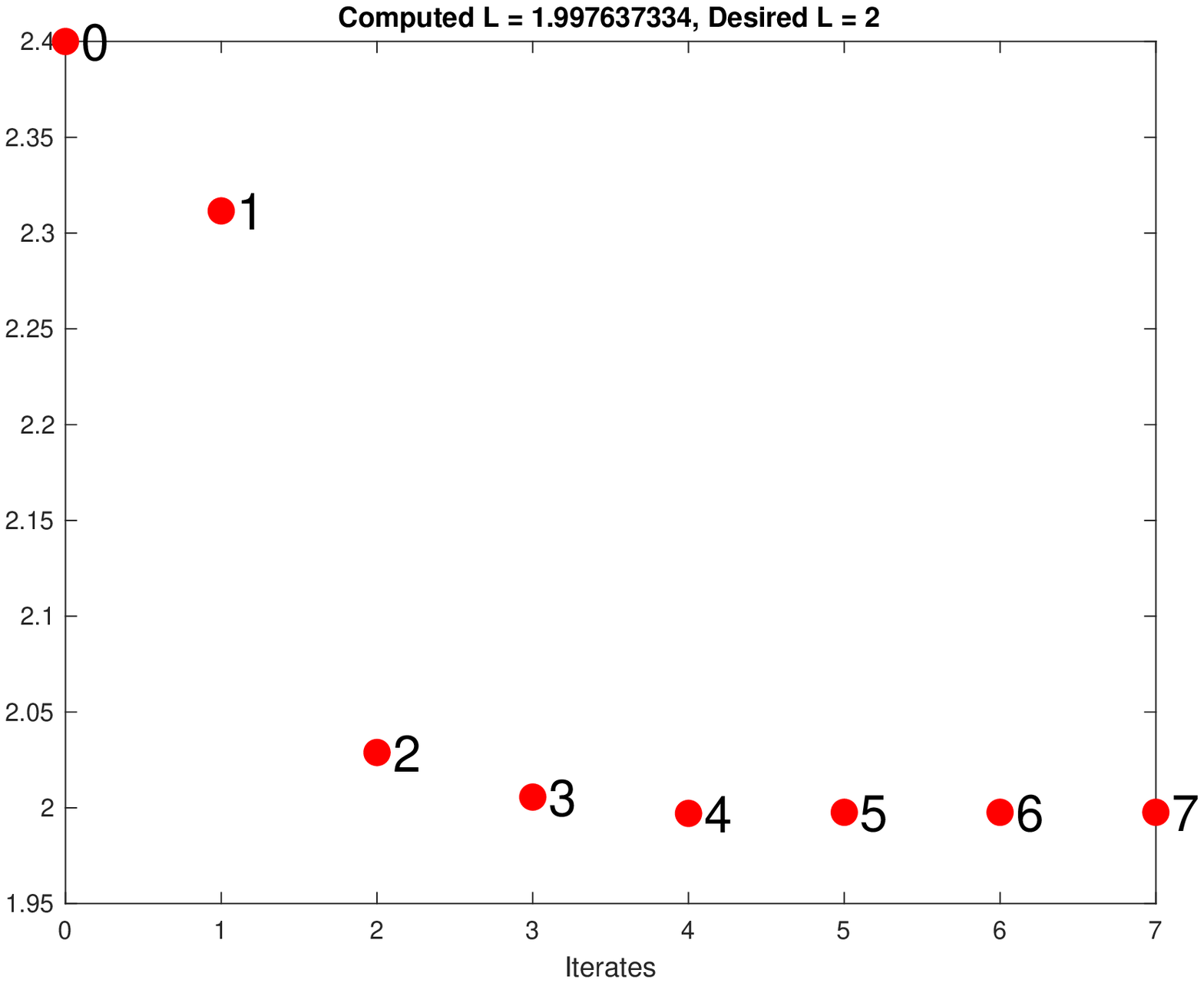}
\caption{Burgers equation, fixed~$u_0$ and~large~$\eta$. The iterates in~\texttt{active-set} algorithm.}
\label{Burgers_iter_case2}
\end{minipage}
\hfill
\begin{minipage}[t]{0.49\linewidth}
\noindent
\includegraphics[width=\linewidth]{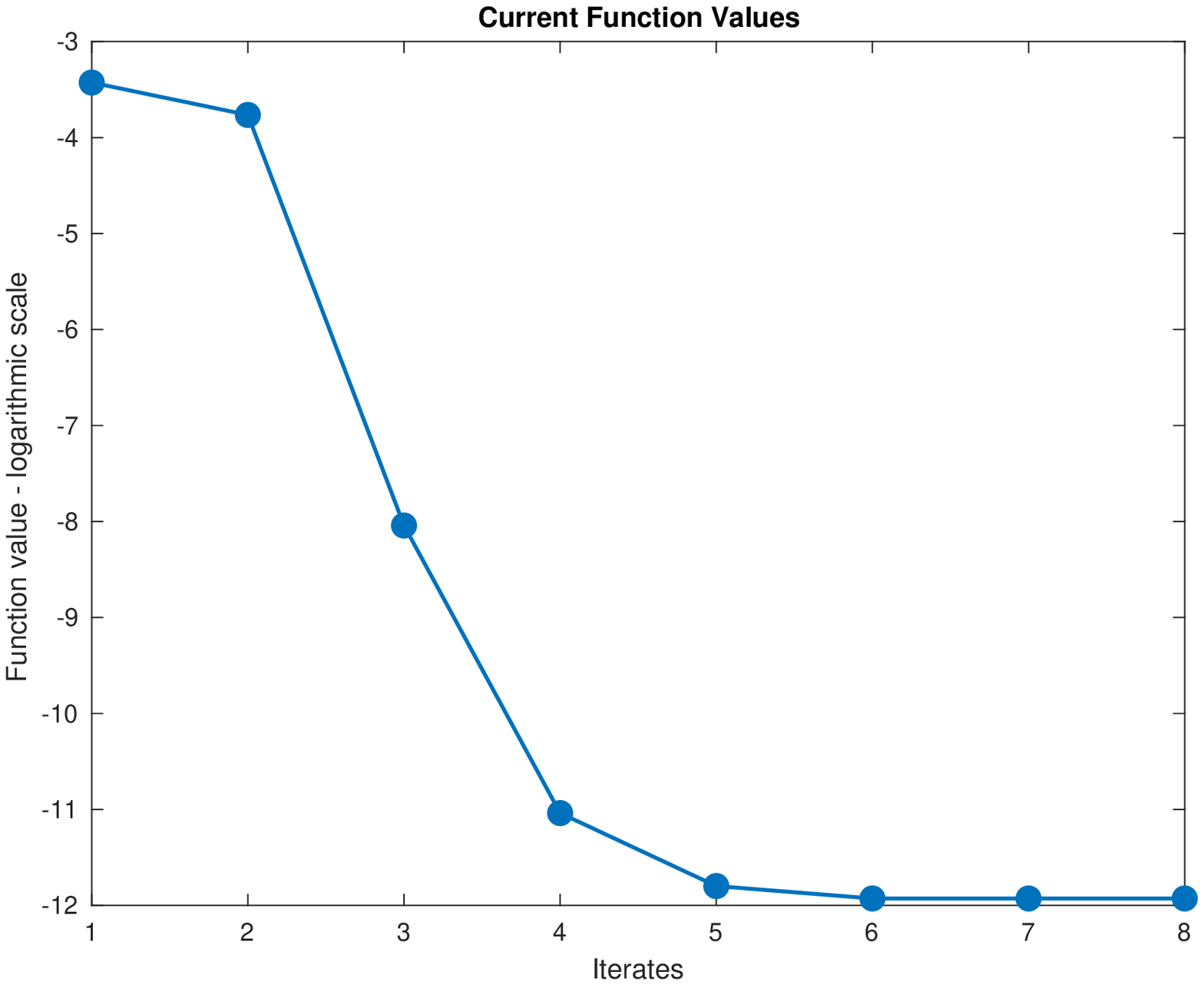}
\caption{Burgers equation, fixed~$u_0$ and~large~$\eta$. Evolution of the cost.}
\label{Burgers_J_case2}
\end{minipage}
\end{figure}

   The results of the numerical implementation are shown in Table~\ref{tab:burgers_case2}, where again random noise was incorporated. The contents of Figures~\ref{Burgers_Sol_case2}, \ref{Burgers_iter_case2} and~\ref{Burgers_J_case2} are similar to those above.
      

%
\begin{table}[h!]
\begin{minipage}[t]{0.46\linewidth}
\centering
\renewcommand{\arraystretch}{1.2}
\medskip
\caption{Burgers equation, fixed $u_0$ and large~$\eta$. Results with random noise in the target
{(desired length: $L_d = 2$). } \\
{ }
{ }}
\begin{tabular}{|cccc|} \hline
\%  noise  & Cost  & Iterates &  Computed $L_c$ \\
\hline
1\%           &   1.e-2         &       6      &   2.032815856  \\ \hline
0.1\%        &    1.e-5        &     11       &  2.012510004\\ \hline
0.01\%      &   1.e-5       &       9          &  1.985859861\\\hline
0.001\%    &   1.e-6       &     9           &   1.994836103\\ \hline
0\%           &   1.e-6      &     9         &    1.997637334  \\  \hline
\end{tabular}
\label{tab:burgers_case2}
\end{minipage}
\hfill
\begin{minipage}[t]{0.54\linewidth}
\vspace{0cm}
\includegraphics[width=\linewidth]{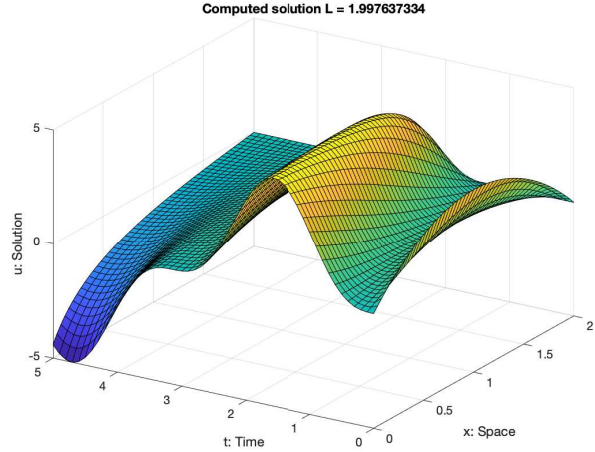}
\captionof{figure}{Burgers equation, $u_0\neq 0$ and~large~$\eta$. The computed solution.}
\label{Burgers_Sol_case2}
\end{minipage}
\end{table}
%

\noindent 
\textbf{Case 1.3: Burgers equation with $u_0\neq 0$ and ``small'' $\eta$.}

   Here, we deal with a non-uniqueness situation.
   Our aim is to investigate the behavior of the algorithm in a situation of this kind.
   
   We take $T=6$, $\eta=0$ in~$(0,T)$ and~$u_0(x) \equiv \pi\sin(\pi x/2)/(2+\cos(\pi x/2))$.
   Note that we have $u_0(x) \equiv \sin(3\pi x/L^1_d)/(2+\cos(3\pi x/L^1_d)) \equiv \sin(2\pi x/L^2_d)/(2+\cos(2\pi x/L^2_d)) $, with $L^1_d = 6$ and~$L^2_d = 4$; consequently, this initial data can be used as in~Section~\ref{SSSec-2.1.2} to prove non-uniqueness.
   
\begin{figure}[h!]
\begin{minipage}[t]{0.47\linewidth}
\noindent
\includegraphics[width=\linewidth]{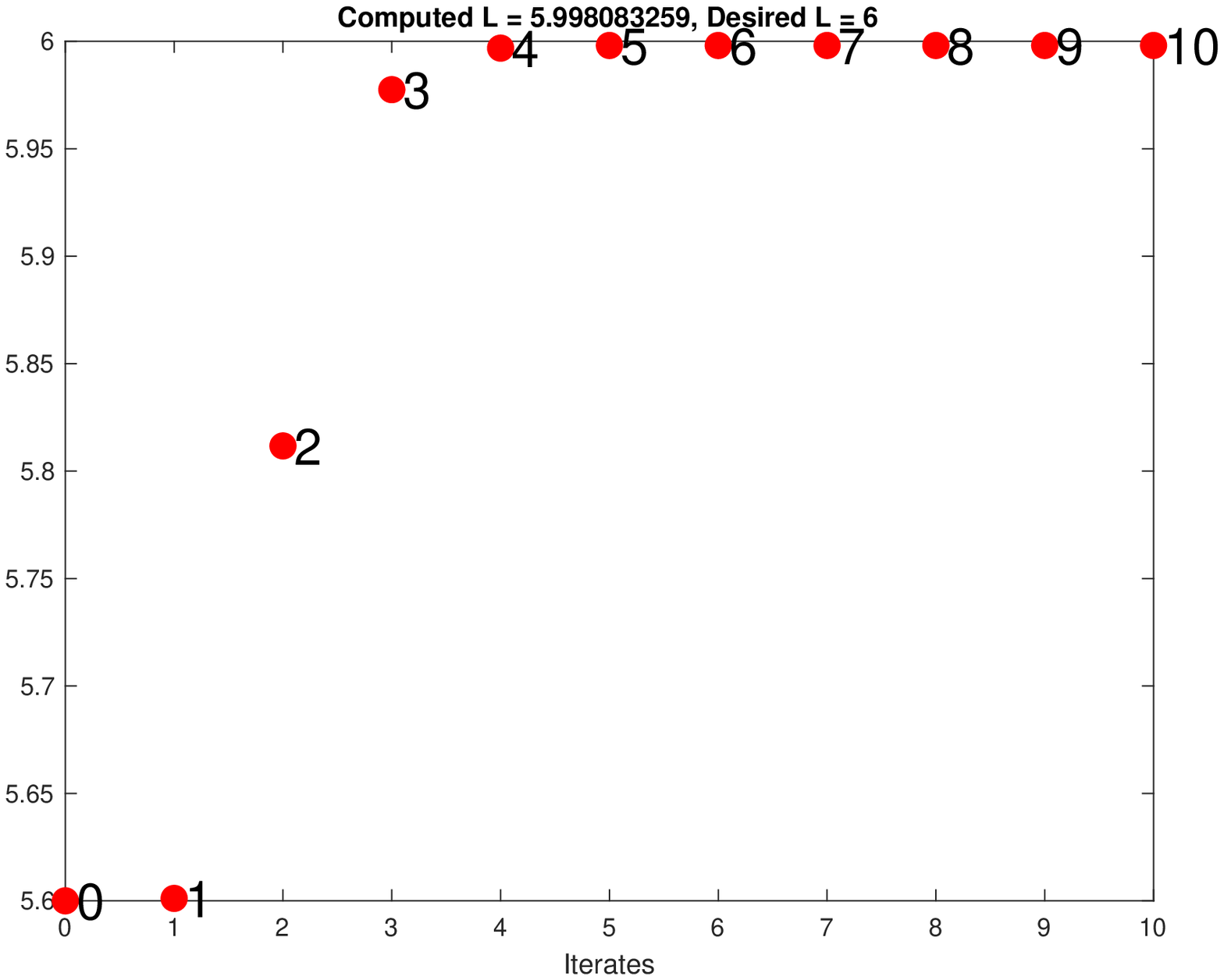}
\caption{Burgers equation, $\eta=0$,  fixed $u_0(x)$. Iterates in \texttt{active-set} algorithm with $L^1_d=6$.}
\label{Burgers_iter_case3}
\end{minipage}
\hfill
\begin{minipage}[t]{0.47\linewidth}
\noindent
\includegraphics[width=\linewidth]{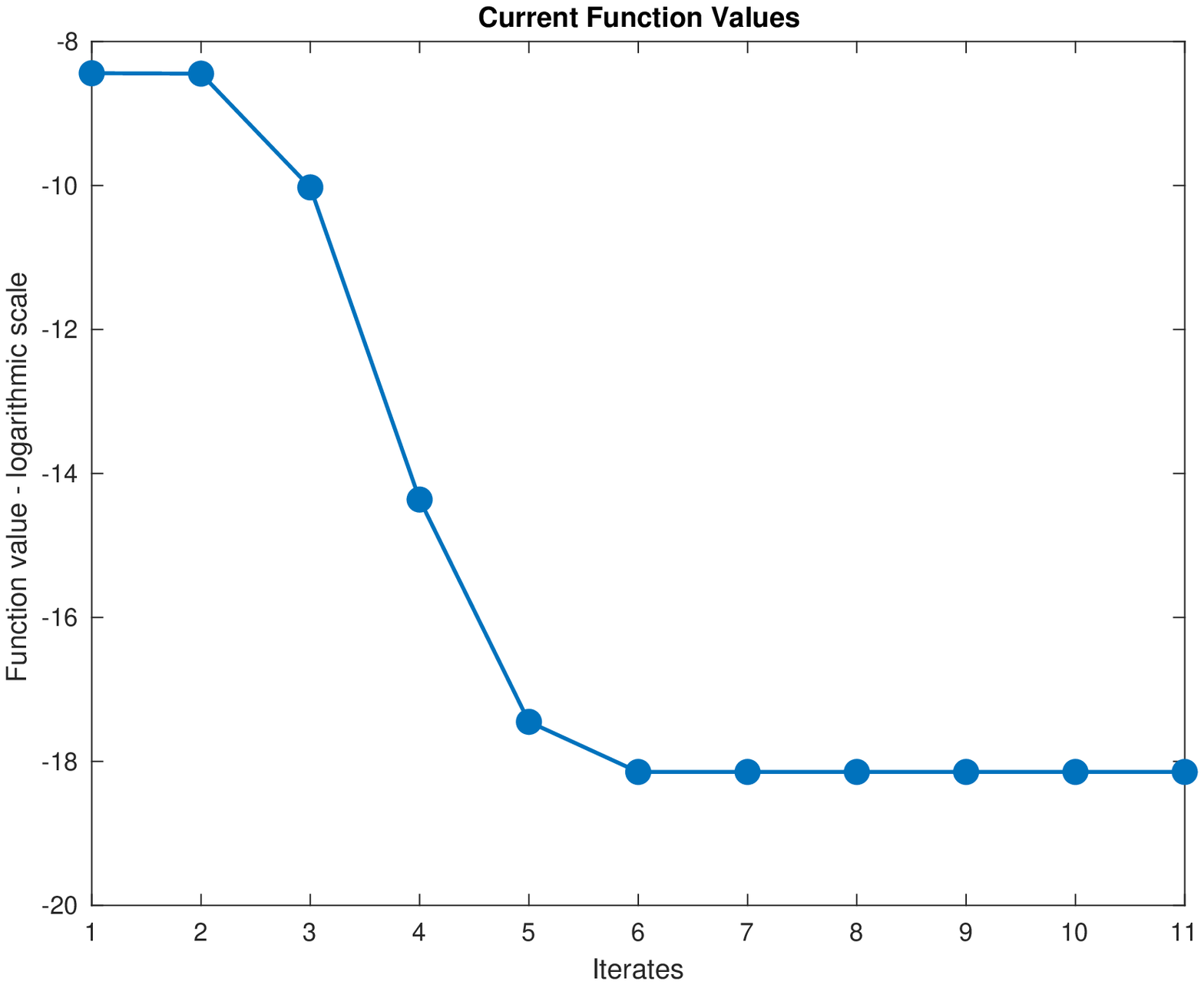}
\caption{Burgers equation, $\eta=0$,  fixed $u_0(x)$. Evolution of the cost for $L^1_d=6$, $J(L^1_c)<10^{-8}$.}
\label{Burgers_J_case3}
\end{minipage}
\end{figure}

\begin{figure}[h!]
\begin{minipage}[t]{0.47\linewidth}
\noindent
\includegraphics[width=\linewidth]{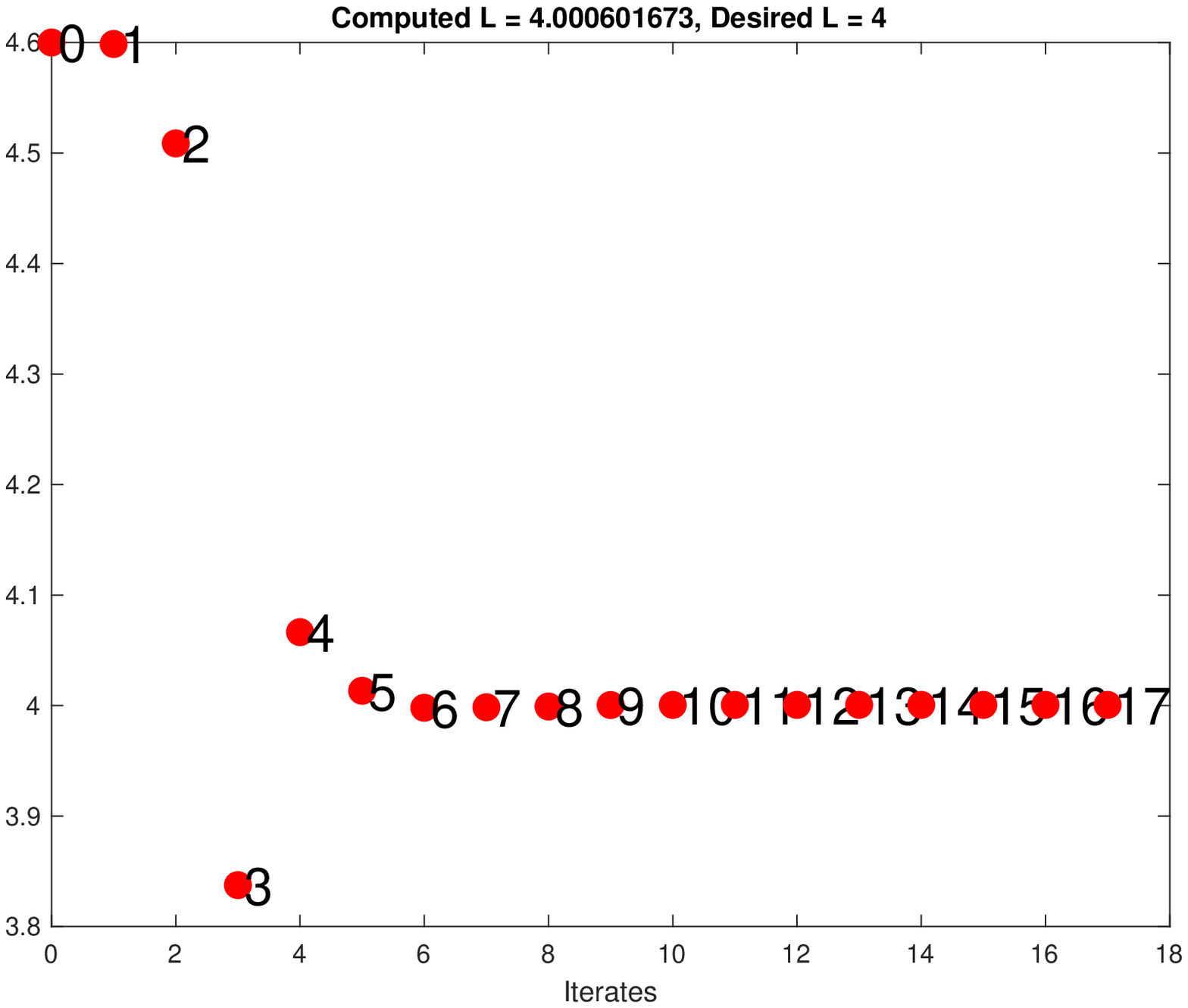}
\caption{Burgers equation, $\eta=0$,  fixed $u_0(x)$. Iterates in \texttt{active-set} algorithm with $L^2_d = 4$.}
\label{Burgers_iter_case3b}
\end{minipage}
\hfill
\begin{minipage}[t]{0.49\linewidth}
\noindent
\includegraphics[width=\linewidth]{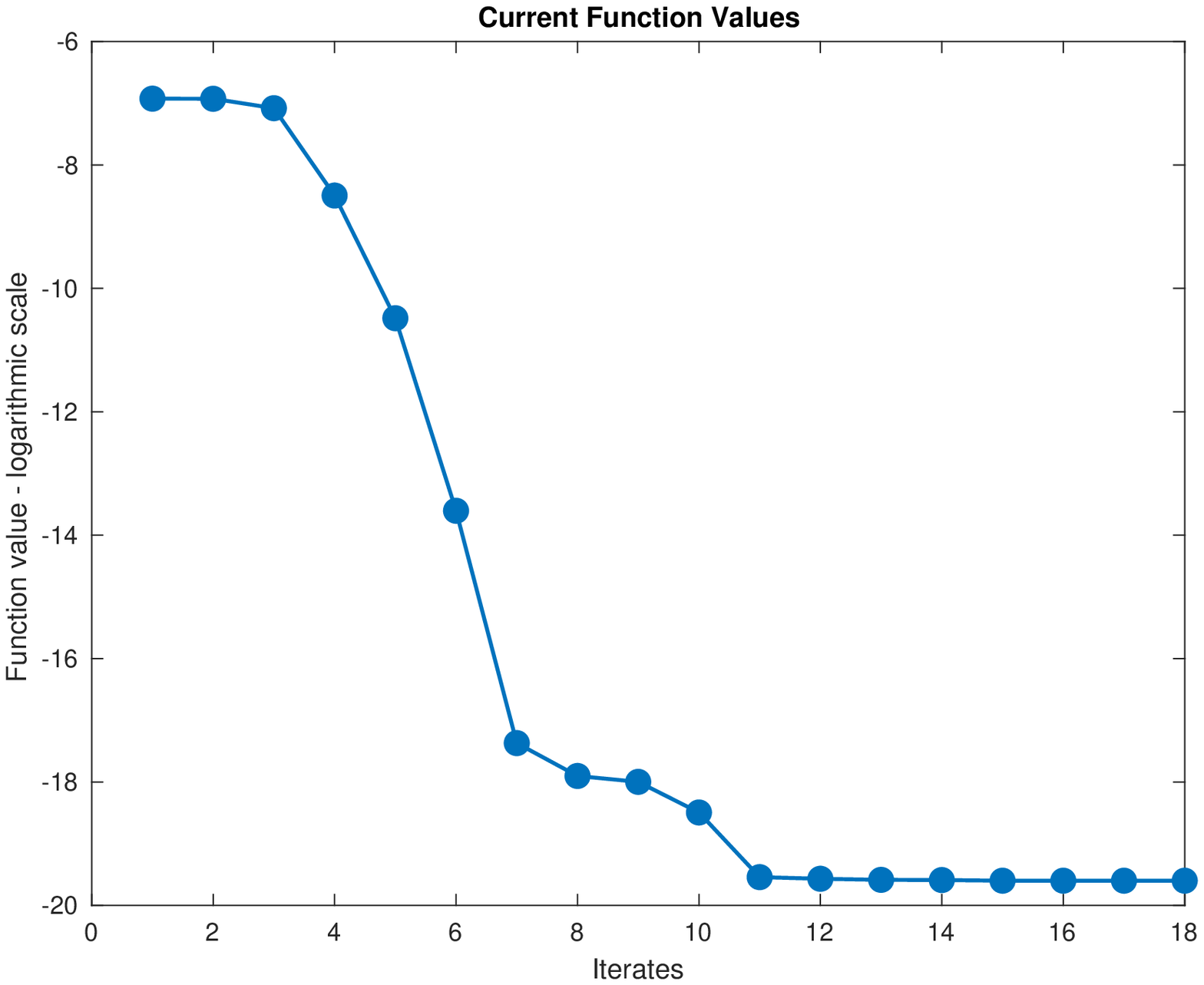}
\caption{Burgers equation, $\eta=0$,  fixed $u_0(x)$. Evolution of the cost for $L^2_d=4$, $J(L_c^2)<10^{-9}$.}
\label{Burgers_J_case3b}
\end{minipage}
\end{figure}

   We will consider the following experiments:
   
\begin{itemize}

\item First, we start {\black from~$L_i = 5.6$,} and we obtain the results exhibited in~Figures~\ref{Burgers_iter_case3} and~\ref{Burgers_J_case3}.
   The computed value is $L^1_c=5.998083259$ and the associated cost is $J(L_c^1)< 10^{-8}$.
\item Then, we start {\black from~$L_i = 4.6$,} and we obtain the results exhibited in~Figures~\ref{Burgers_iter_case3b} and~\ref{Burgers_J_case3b}.
   The computed value is $L^2_c=4.000601673$ and the associated cost is again $J(L_c^2)< 10^{-9}$.
   
\end{itemize}

\begin{figure}[h!]
\begin{minipage}[t]{0.47\linewidth}
\noindent
\includegraphics[width=\linewidth]{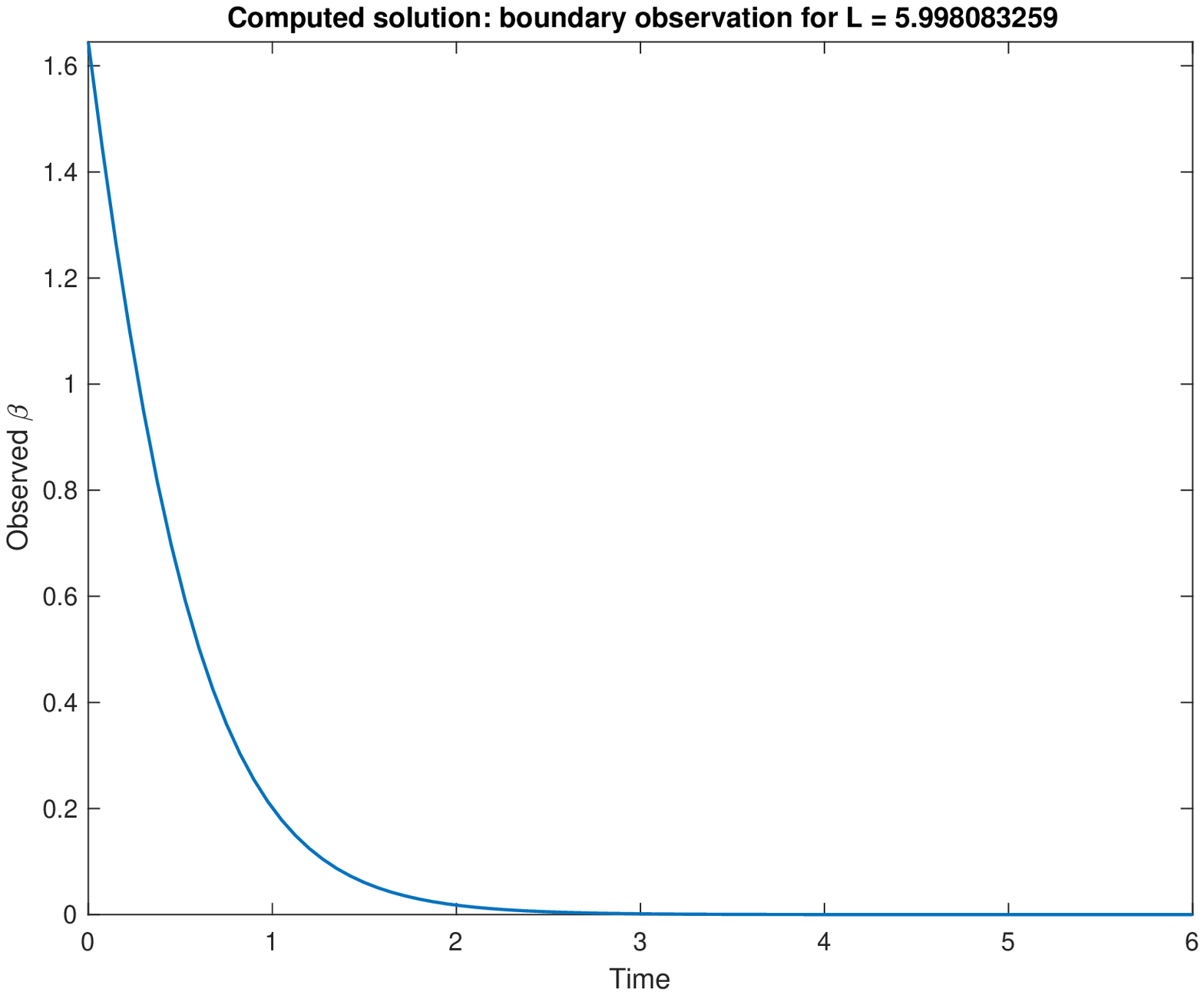}
\caption{Burgers equation, $\eta=0$,  fixed $u_0(x)$. The computed boundary observation $u_x(0,\cdot)$ for~$L^1_c=5.996562049$.}
\label{Burgers_Sol_case3}
\end{minipage}
\hfill
\begin{minipage}[t]{0.47\linewidth}
\noindent
\includegraphics[width=\linewidth]{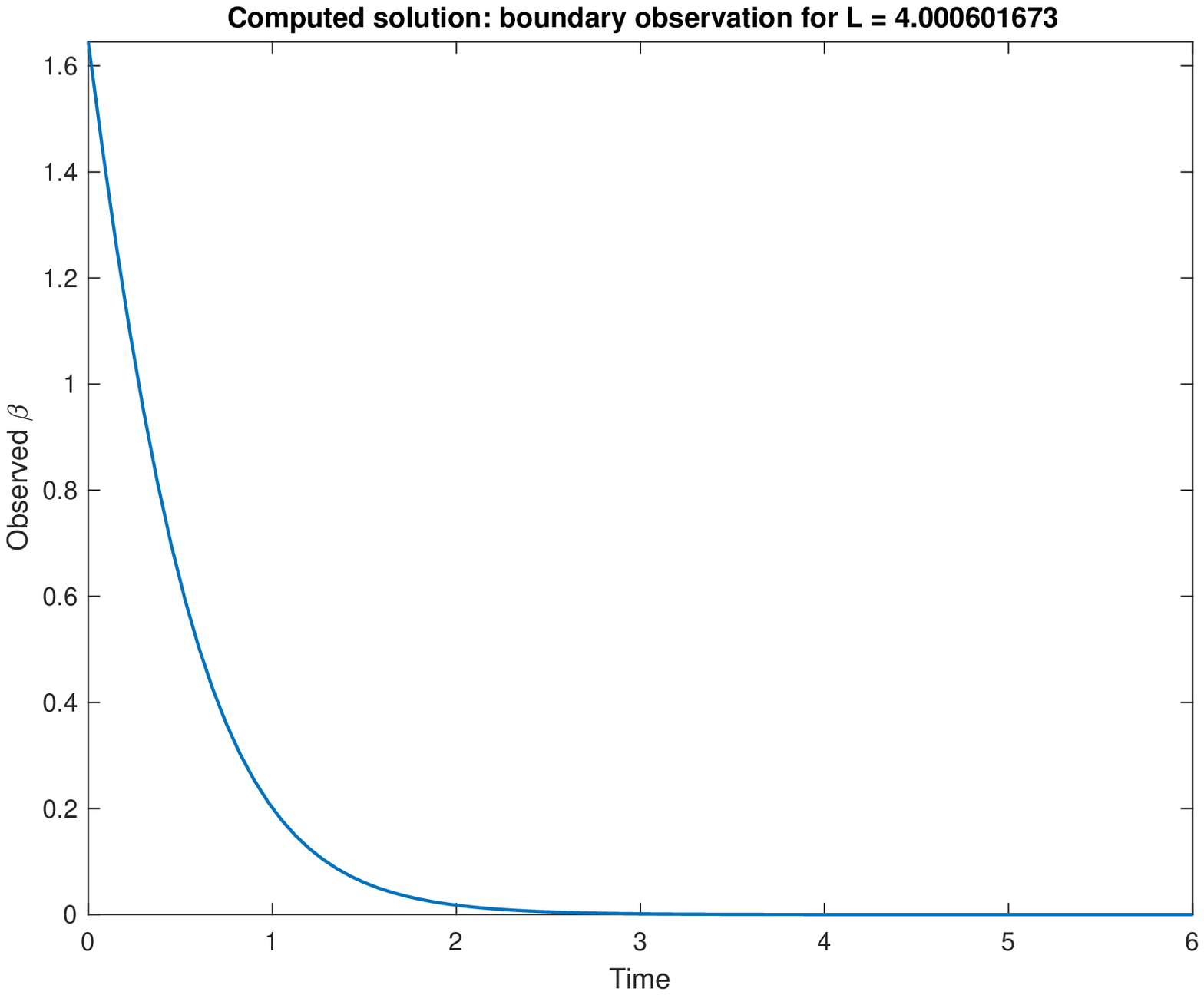}
\caption{Burgers equation, $\eta=0$,  fixed $u_0(x)$. The computed boundary observation $u_x(0,\cdot)$ for~$L_c^2=4.007345905$}
\label{Burgers_Sol_case3b}
\end{minipage}
\end{figure}

   The corresponding computed boundary observations are displayed in~Figures~\ref{Burgers_Sol_case3} and~\ref{Burgers_Sol_case3b}, respectively.
   Thus, we confirm that these identical observations correspond, as we already knew, two different solutions. 
 
\begin{figure}[h!]
\begin{minipage}[t]{0.49\linewidth}
\noindent
\includegraphics[width=\linewidth]{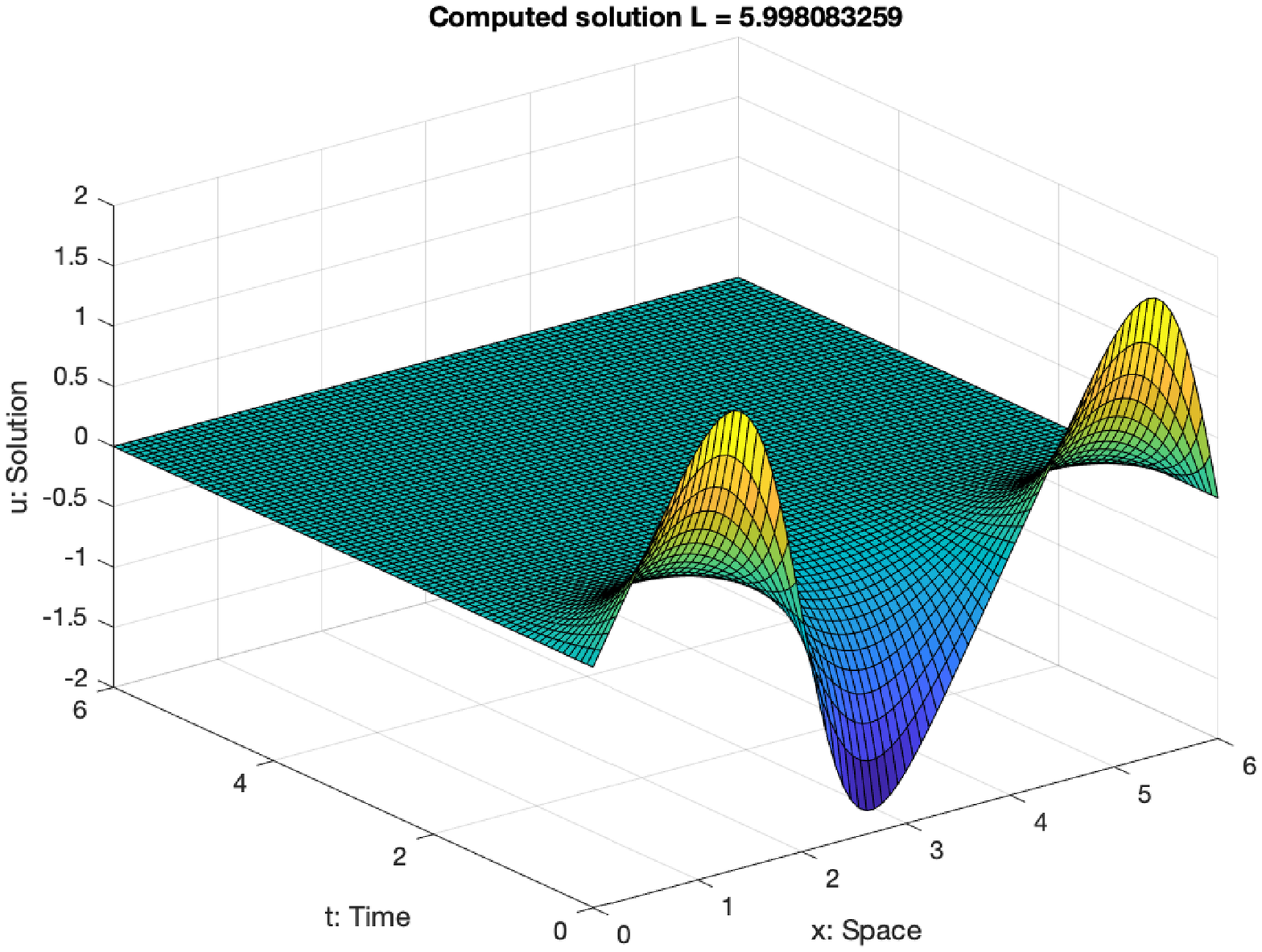}
\caption{Burgers equation, $\eta=0$,  fixed $u_0(x)$. The computed solution corresponding to $L^1_c=5.998083259$.}
\label{Burgers_Sol_case3Bis}
\end{minipage}
\hfill
\begin{minipage}[t]{0.49\linewidth}
\noindent
\includegraphics[width=\linewidth]{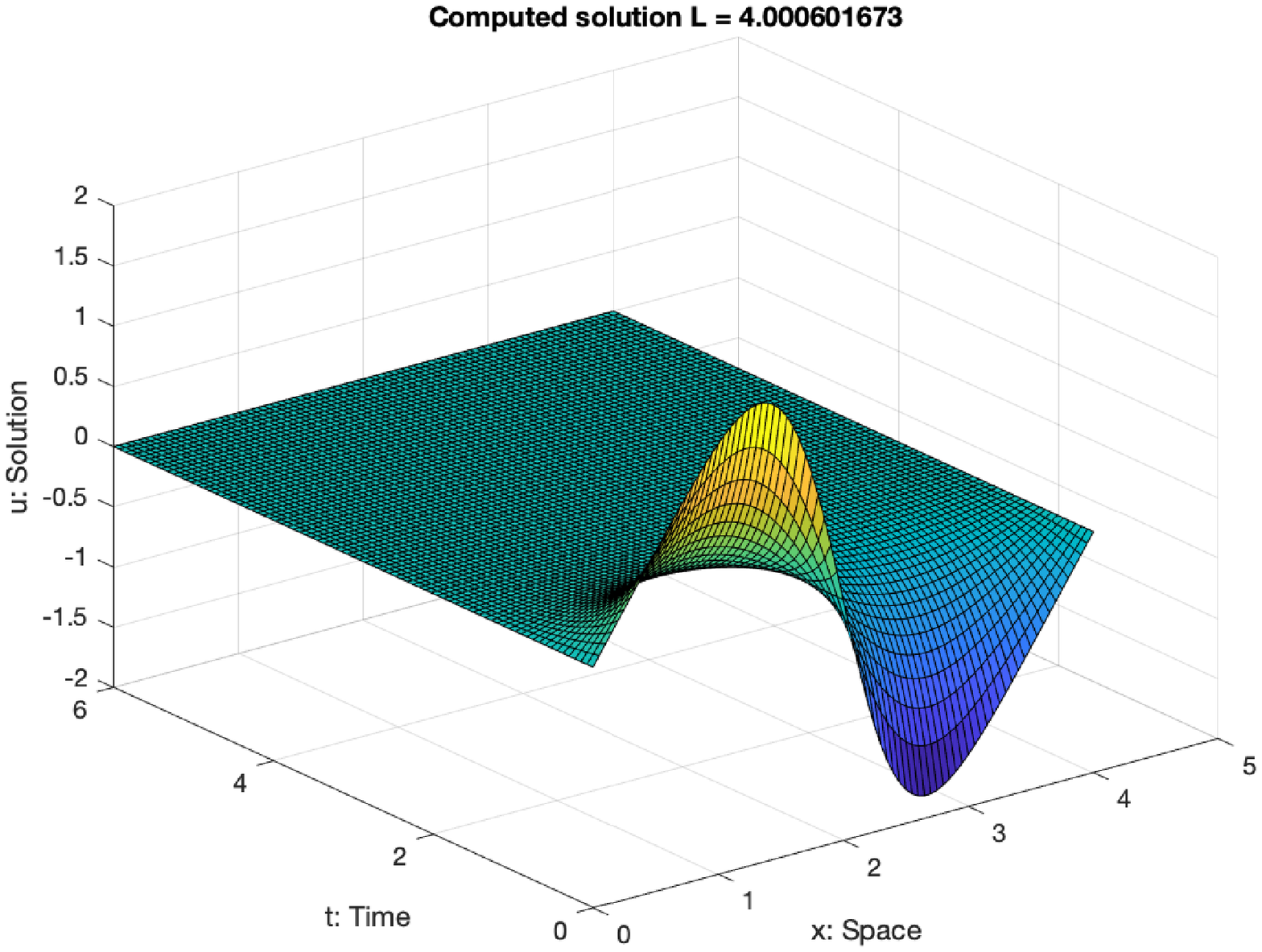}
\caption{Burgers equation, $\eta=0$,  fixed $u_0(x)$. The computed solution corresponding to  $L_c^2=4.000601673$.}
\label{Burgers_Sol_case3bBis}
\end{minipage}
\end{figure}

\subsection{Inverse problems for the Burgers-heat system}

   This section is concerned with \textbf{IP-2} and other related problems.
   We will consider several choices of boundary conditions and also several different observations.


\subsubsection{Dirichlet boundary conditions for $u$ and $\theta$ and stress and flux observations.}

   We consider the system~\eqref{pbBH}.
   A reformulation of~\textbf{IP-2} is the following:
   \[
\left\{
\begin{array}{l}\dis
\text{Minimize } J_2(\ell) := \frac{1}{2} \int_0^T |\beta(t) - u_x^\ell(0,t)|^2 \,dt 
+ \frac{1}{2} \int_0^T |\alpha(t) - \theta_x^\ell(0,t)|^2 \,dt
\\ \noalign{\smallskip} \dis
\text{Subject to: } \ell \in (\ell_0,\ell_1), \ \text{ $(u^\ell,\theta^\ell)$ satisfies~\eqref{pbBH}.}
\end{array}
\right.
   \]
\noindent
\textbf{Case 2.1: Burgers-heat system with $(u_0,\theta_0)=(0,0)$ and $\eta\neq 0$ and $\lambda\neq0$.}

   We take $T=5, \ \ \eta(t) \equiv 5\sin^3 t, \ \ \lambda(t) \equiv 0.2\cos(t)\sin(t) \ \text{ and } \ (u_0(x),\theta_0(x))\equiv (0,0).
$
   Starting from $L_{i} = 1$, our goal is to recover the desired value of the length $L_d= 2$.

\begin{figure}[h!]
\begin{minipage}[t]{0.49\linewidth}
\vspace{0cm}
\includegraphics[width=\linewidth]{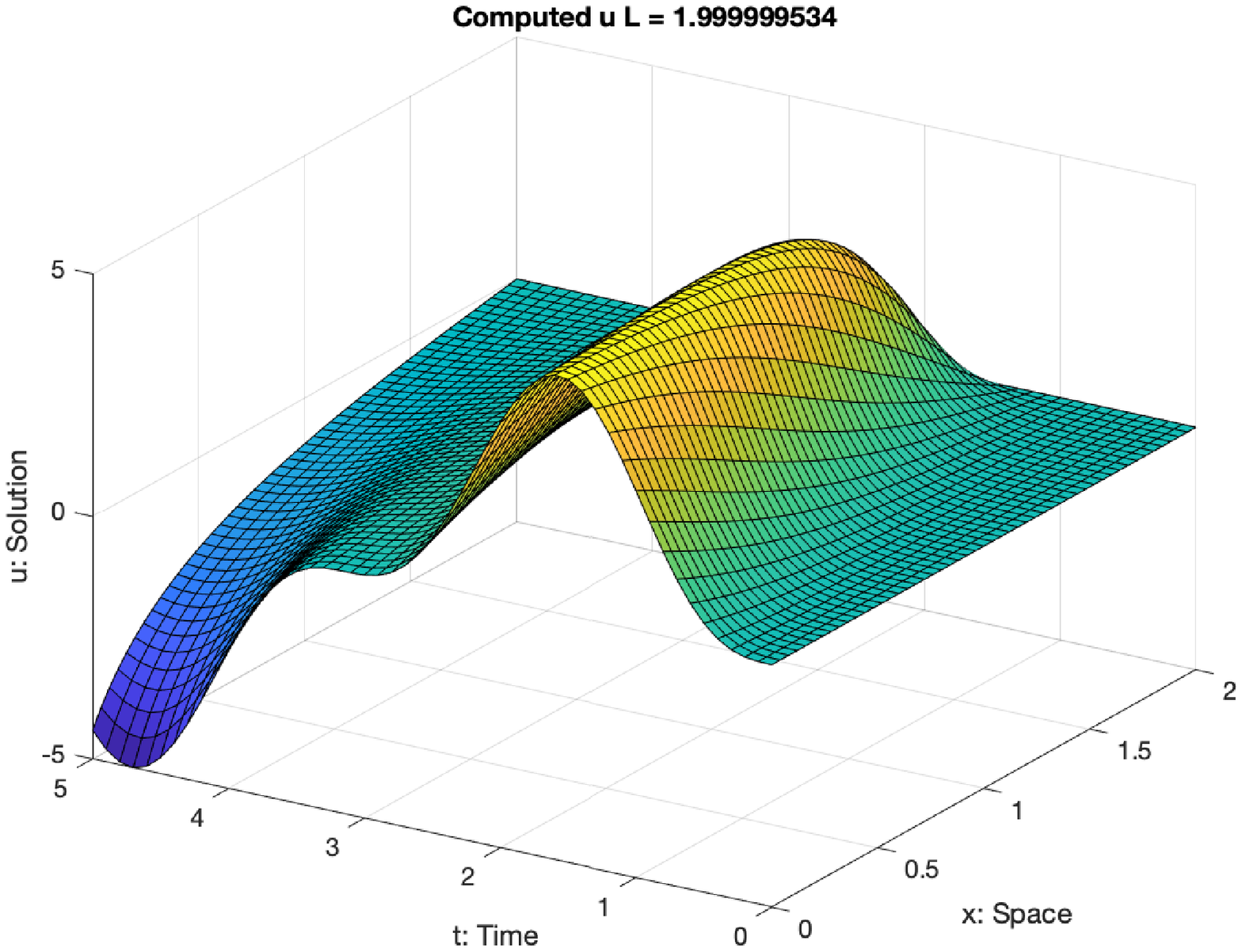}
\captionof{figure}{Burgers equation with heat effect with $(u_0,\theta_0)=(0,0)$ and $(\eta,\lambda)\neq (0,0)$ with two observations $u_x(0,t)$ and $\theta_x(0,t)$. The computed solution $u$.}
\label{BH_DDSoluC5}
\end{minipage}
\hfill
\begin{minipage}[t]{0.49\linewidth}
\vspace{0cm}
\includegraphics[width=\linewidth]{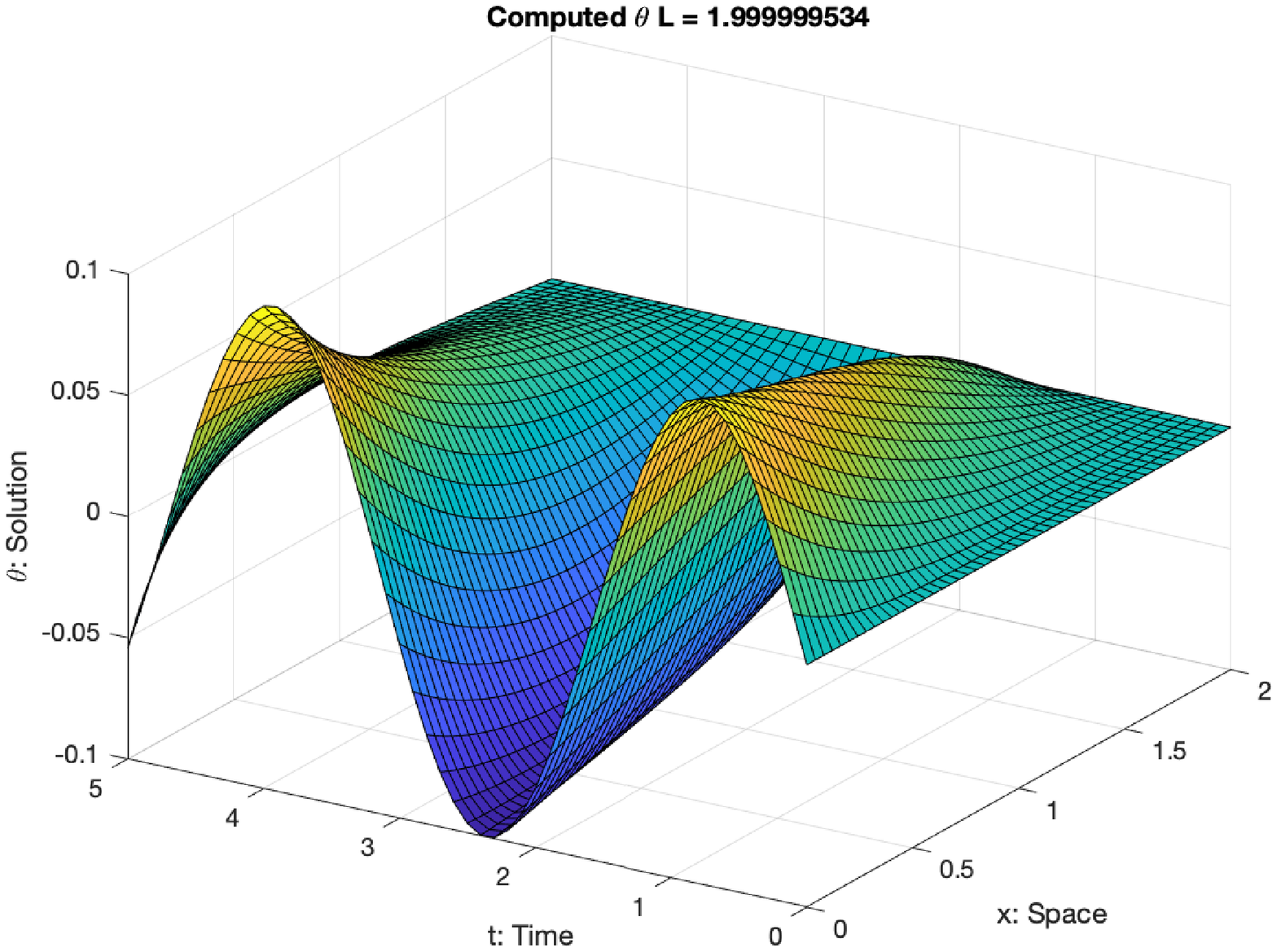}
\captionof{figure}{Burgers equation with heat effect with $(u_0,\theta_0)=(0,0)$ and$(\eta,\lambda)\neq (0,0)$ with two observations $u_x(0,t)$ and $\theta_x(0,t)$. The computed solution $\theta$.}
\label{BH_DDSolthC5}
\end{minipage}
\end{figure}

\begin{figure}[h!]
\begin{minipage}[t]{0.49\linewidth}
\noindent
\includegraphics[width=\linewidth]{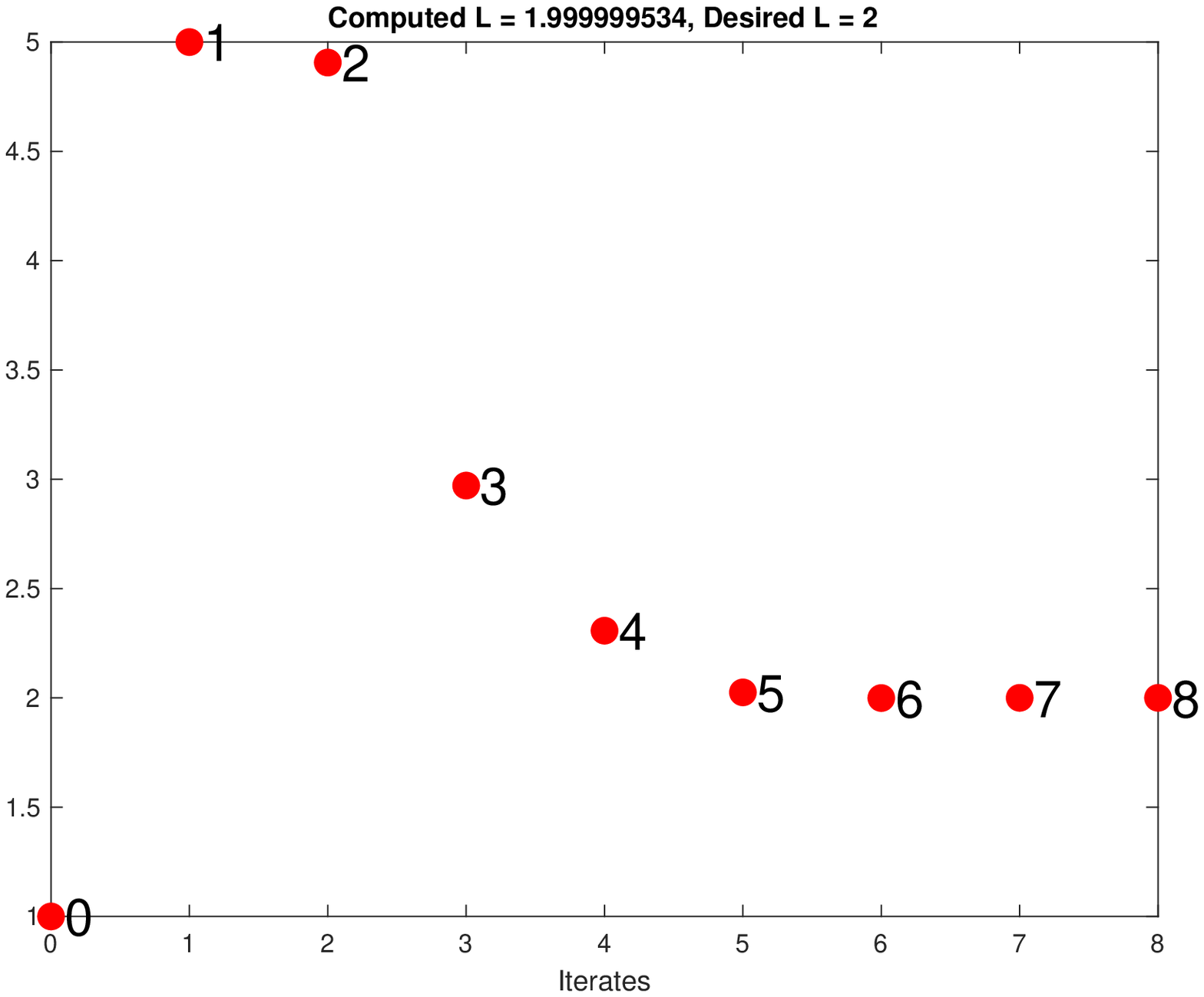}
\caption{Burgers equation with heat effect with $(u_0,\theta_0)=(0,0)$ and $(\eta,\lambda)\neq (0,0)$ with two observations $u_x(0,t)$ and $\theta_x(0,t)$. The iterates in \texttt{active-set} algorithm.}
\label{BH_DDiterC5}
\end{minipage}
\hfill
\begin{minipage}[t]{0.49\linewidth}
\noindent
\includegraphics[width=\linewidth]{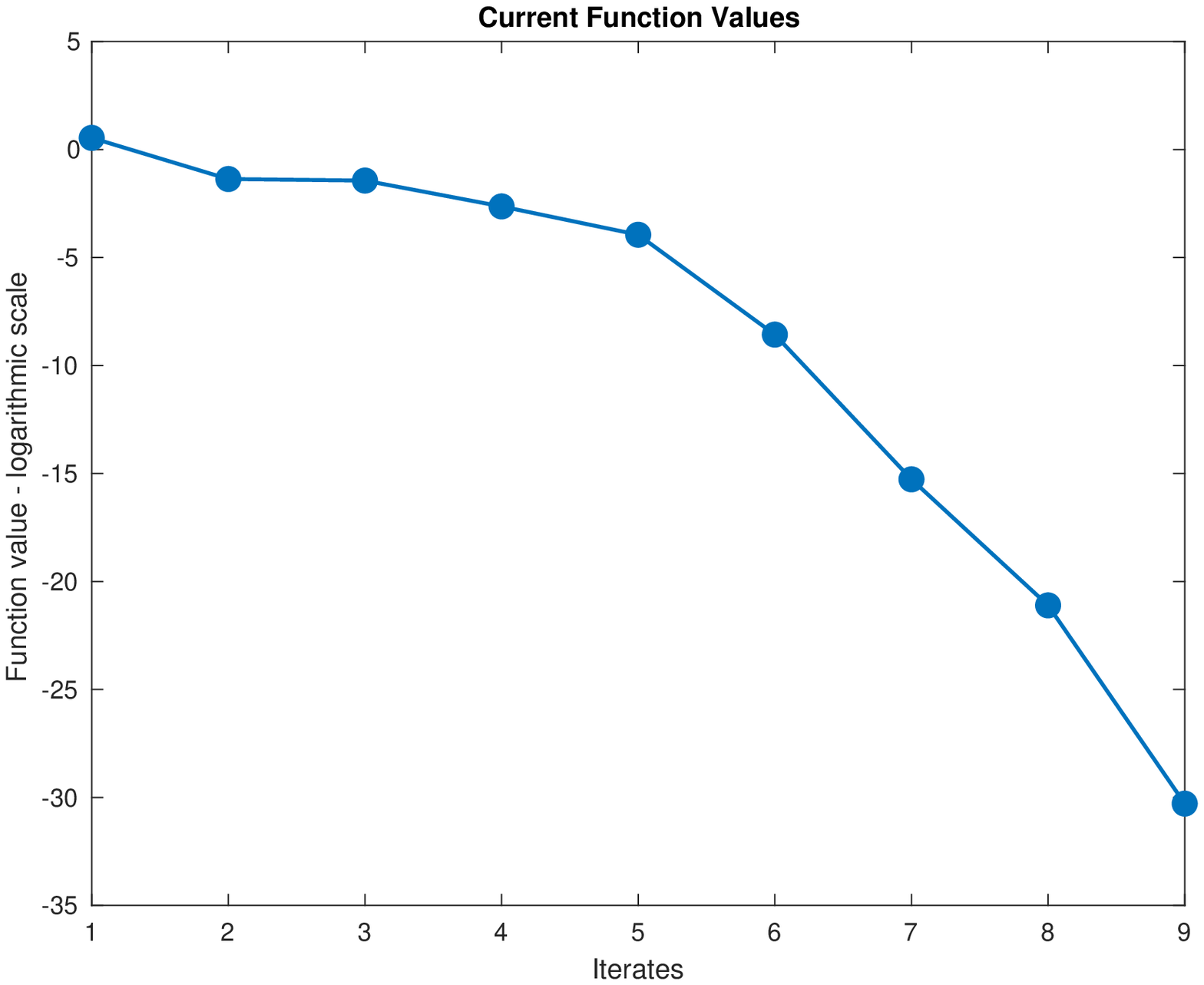}
\caption{Burgers equation with heat effect with $(u_0,\theta_0)=(0,0)$ and $(\eta,\lambda)\neq (0,0)$ with two observations $u_x(0,t)$ and $\theta_x(0,t)$. Evolution of the cost.}
\label{BH_DDfvC5}
\end{minipage}
\end{figure}

   The computed length is~$L_c=1.999999534$, the cost is $J(L_c)< 10^{-14}$ is reached at the iterate 8 of the optimization algorithm.
   The corresponding solution to~\eqref{pbBH} is displayed in~Figures~\ref{BH_DDSoluC5} and~\ref{BH_DDSolthC5}.
   The evolution of the iterates and the cost in the minimization process in the absence of the random noise appear in~Figures~\ref{BH_DDiterC5} and~\ref{BH_DDfvC5}, respectively.


%
%

\noindent 
\textbf{Case 2.2: Burgers-heat system with $(u_0,\theta_0)\neq (0,0)$ and large $\eta$.}

   We take $T=5$,  $\eta(t) = 5\sin^3 t$ and $\lambda(t) = 6\sin(t)\cos(t)$ in~$(0,T)$, $u_0(x) \equiv 0.1 x(2-x)$ and~$\theta_0(x) \equiv 0.1x^2(x-3)$.
   Starting from $L_{i} = 1$, our goal is to recover the desired value of the length $L_d= 2$.

\begin{figure}[h!]
\begin{minipage}[t]{0.49\linewidth}
\vspace{0cm}
\includegraphics[width=\linewidth]{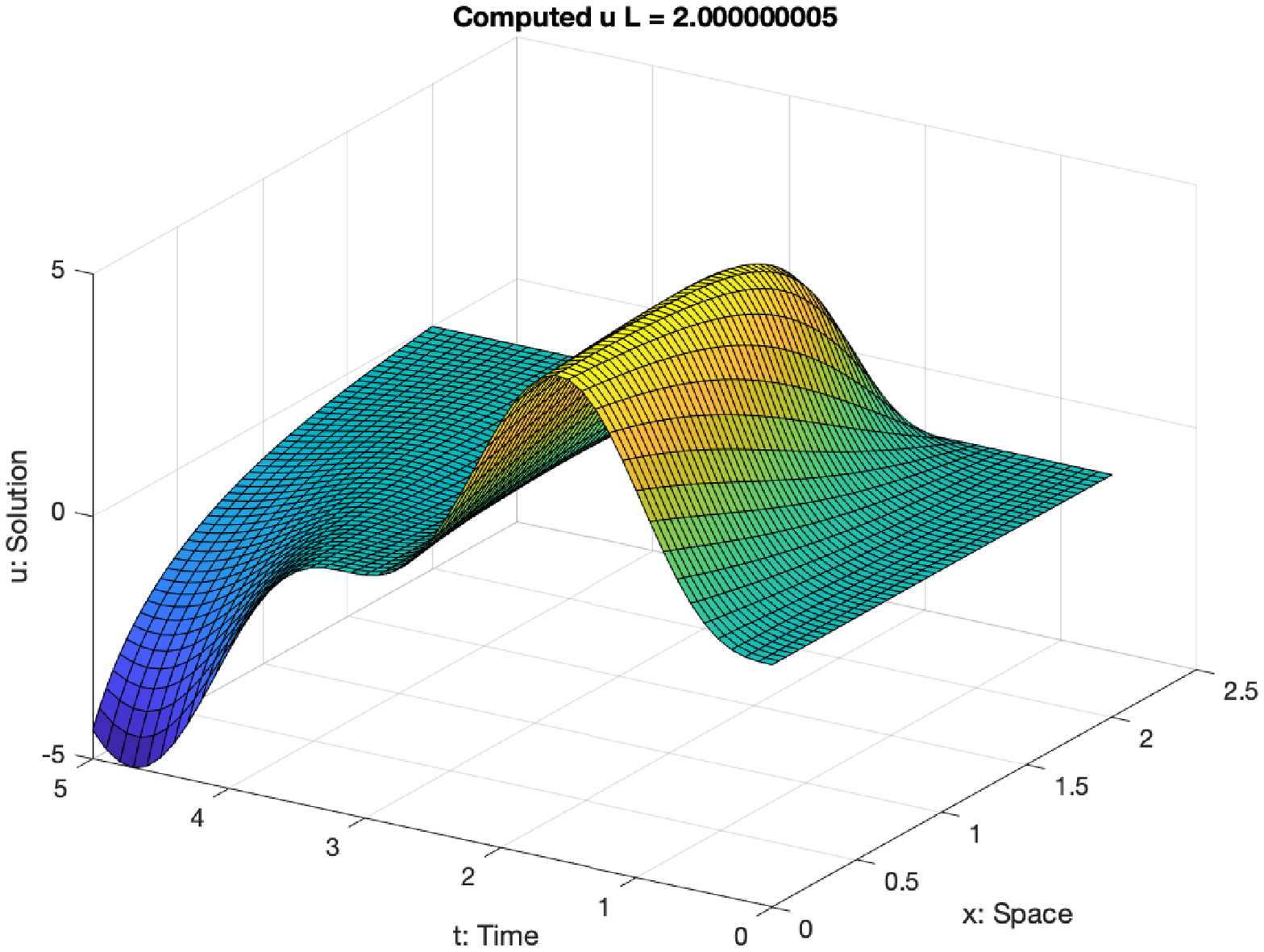}
\captionof{figure}{Burgers equation with heat effect with $(u_0,\theta_0)\neq (0,0)$ and large  $(\eta,\lambda)$ with two observations $u_x(0,t)$ and $\theta_x(0,t)$. The computed solution $u$.}
\label{BH_DDSoluC6}
\end{minipage}
\hfill
\begin{minipage}[t]{0.49\linewidth}
\vspace{0cm}
\includegraphics[width=\linewidth]{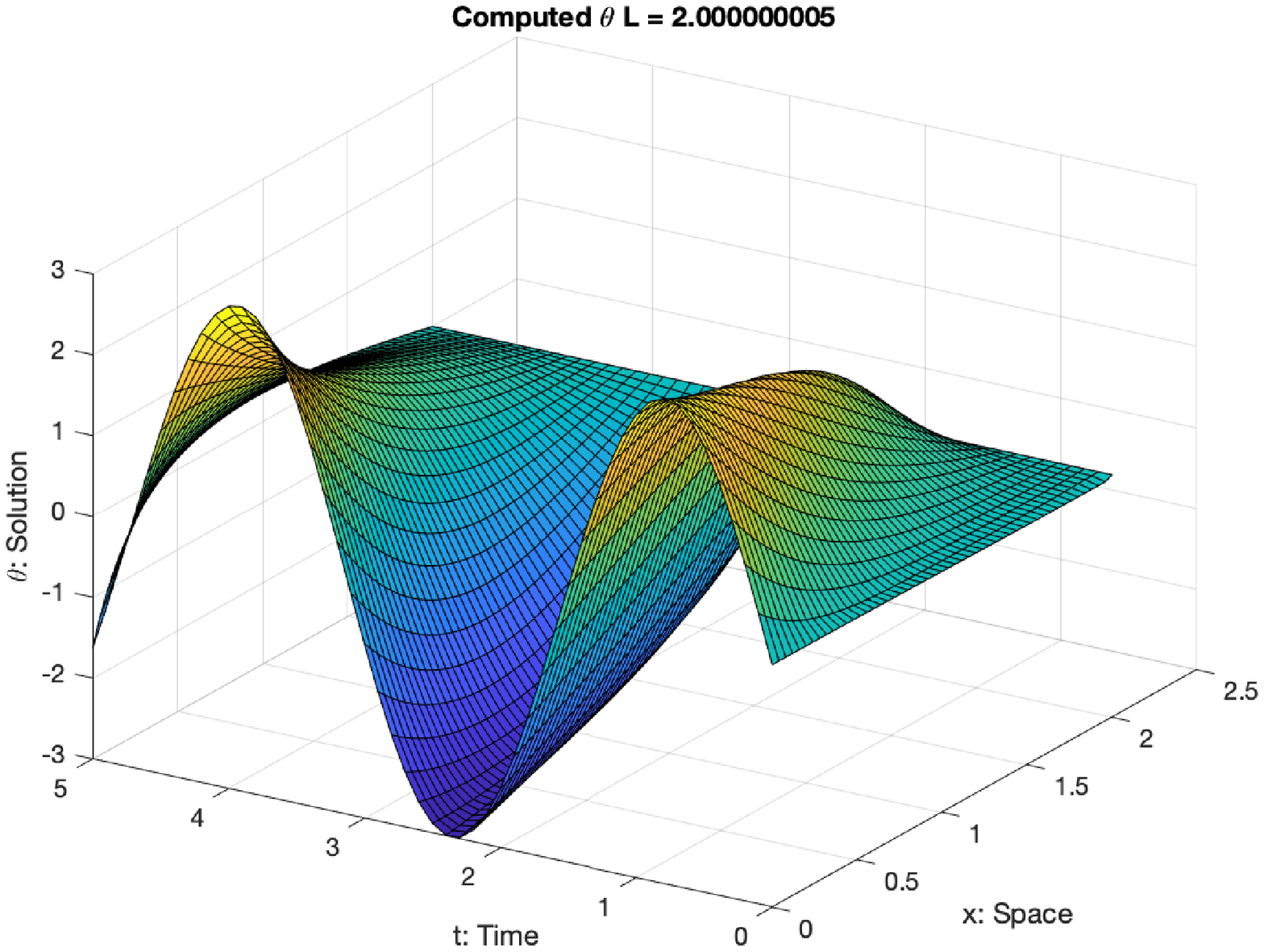}
\captionof{figure}{Burgers equation with heat effect with $(u_0,\theta_0)\neq (0,0)$ and large  $(\eta,\lambda)$ with two observations $u_x(0,t)$ and $\theta_x(0,t)$. The computed solution $\theta$.}
\label{BH_DDSolthC6}
\end{minipage}
\end{figure}

\begin{figure}[h!]
\begin{minipage}[t]{0.49\linewidth}
\noindent
\includegraphics[width=\linewidth]{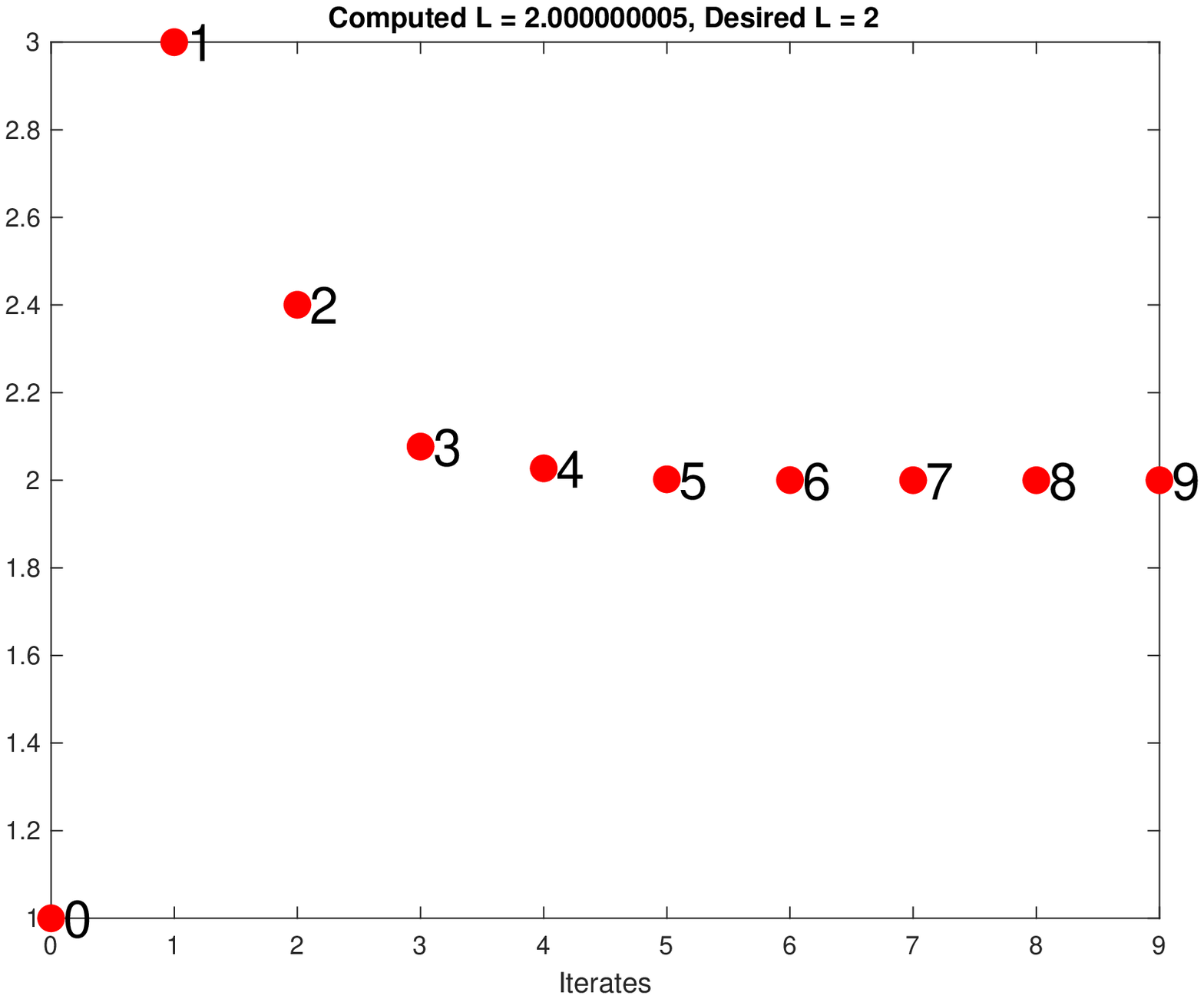}
\caption{Burgers equation with heat effect with $(u_0,\theta_0)\neq (0,0)$ and large  $(\eta,\lambda)$ with two observations $u_x(0,t)$ and $\theta_x(0,t)$. The iterates in \texttt{active-set} algorithm.}
\label{BH_DDiterC6}
\end{minipage}
\hfill
\begin{minipage}[t]{0.49\linewidth}
\noindent
\includegraphics[width=\linewidth]{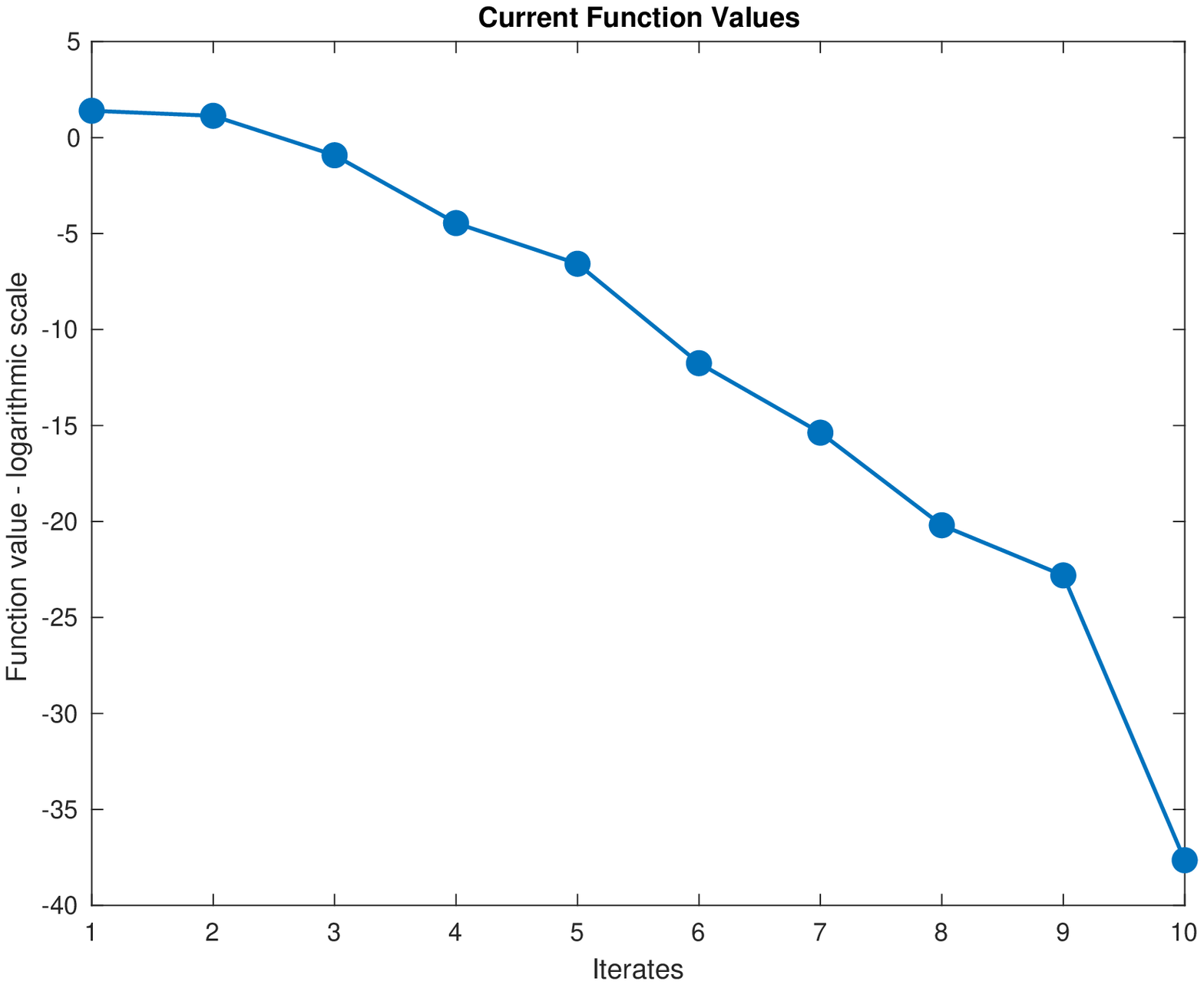}
\caption{Burgers equation with heat effect with $(u_0,\theta_0)\neq (0,0)$ and large  $(\eta,\lambda)$ with two observations $u_x(0,t)$ and $\theta_x(0,t)$. Evolution of the cost.}
\label{BH_DDfvC6}
\end{minipage}
\end{figure}

   The computed length is~$L_c=2.000000005$, the cost is $J(L_c)< 10^{-17}$ is reached at the iterate 9 of the optimization algorithm.
   The corresponding solution to~\eqref{pbBH} is displayed in~Figures~\ref{BH_DDSoluC6} and~\ref{BH_DDSolthC6}.  The evolution of the iterates and the cost in the minimization process in the absence of the random noise appear in~Figures~\ref{BH_DDiterC6} and~\ref{BH_DDfvC6}, respectively.
   

%
%

\subsubsection{Dirichlet boundary conditions for $u$, Neumann boundary conditions for $\theta$ and stress observation.}

   In this section, the system under study is~\eqref{pb-2}.
   The inverse problem is similar to~\textbf{IP-2} and a suitable reformulation is:
   \[
\left\{
\begin{array}{l}\dis
\text{Minimize } J_3(\ell) := \frac{1}{2} \int_0^T |\beta(t) - u_x^\ell(0,t)|^2 \,dt 
\\ \noalign{\smallskip} \dis
\text{Subject to: } \ell \in (\ell_0,\ell_1), \ \text{ $(u^\ell,\theta^\ell)$ satisfies~\eqref{pb-2}.}
\end{array}
\right.
   \]
   As before, two different situations will be analyzed for this problem.
   In both cases, respectively corresponding to zero initial data and nonzero initial data and sufficiently large $\eta$, we will check that uniqueness holds.

\

\noindent
\textbf{Case 2.3: Burgers-heat system with $(u_0,\theta_0)=(0,0)$ and $\eta\neq 0$.}

We observe that this case is reduced to the Burgers single equation. 

 \noindent 
 We take $T=5$,  $\eta(t) = 5\sin^3 t$ in~$(0,T)$ and~$(u_0(x),\theta_0(x))\equiv (0,0)$.
   Starting from $L_{i} = 1$, our goal is to recover the desired value of the length $L_d= 2$.

\begin{figure}[h!]
\begin{minipage}[t]{0.49\linewidth}
\vspace{0cm}
\includegraphics[width=\linewidth]{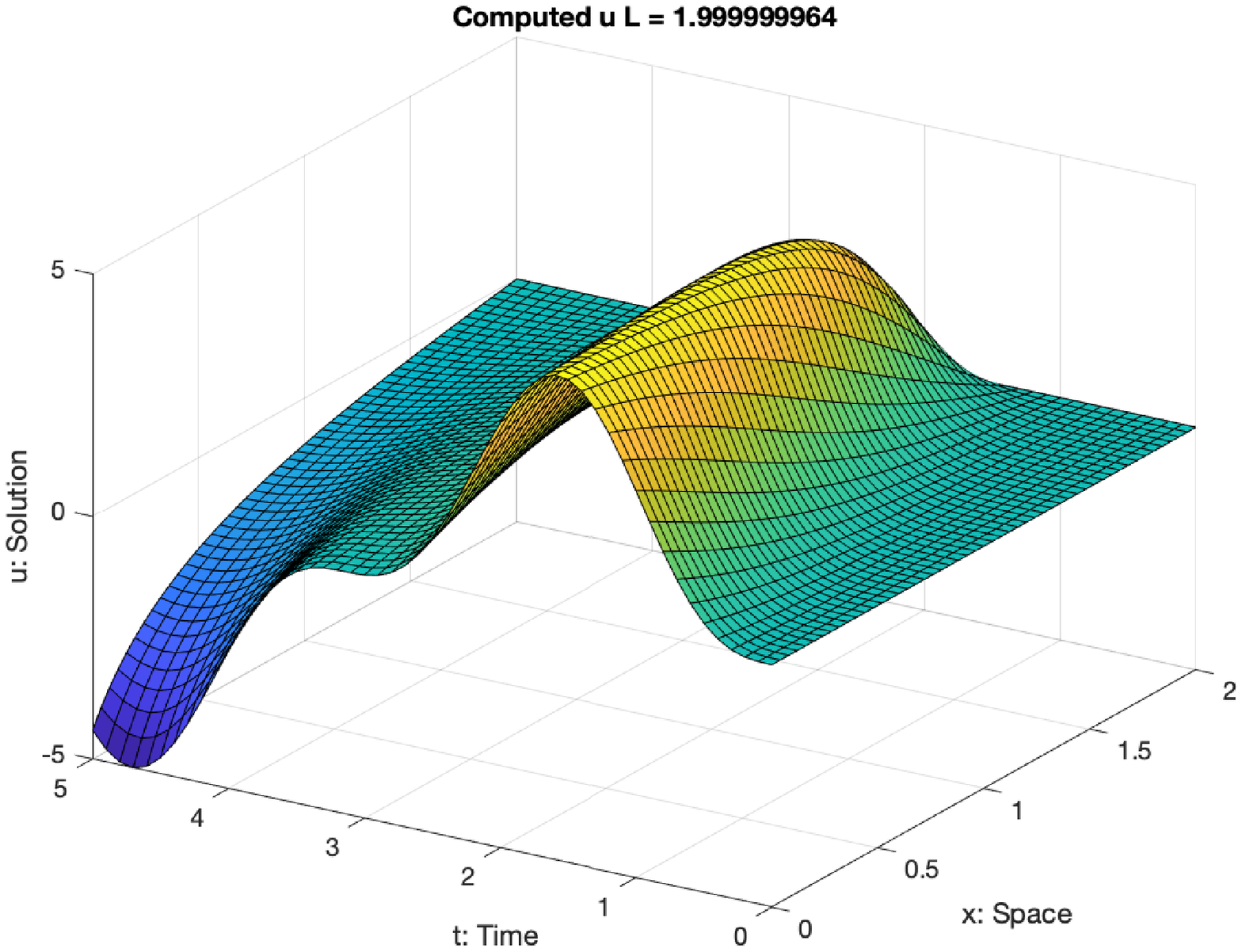}
\captionof{figure}{Burgers equation with heat effect with $(u_0,\theta_0)=(0,0)$ and $\eta\neq 0$ with one observation $u_x(0,t)$. The computed solution $u$.}
\label{BH_DNSoluC1}
\end{minipage}
\hfill
\begin{minipage}[t]{0.49\linewidth}
\vspace{0cm}
\includegraphics[width=\linewidth]{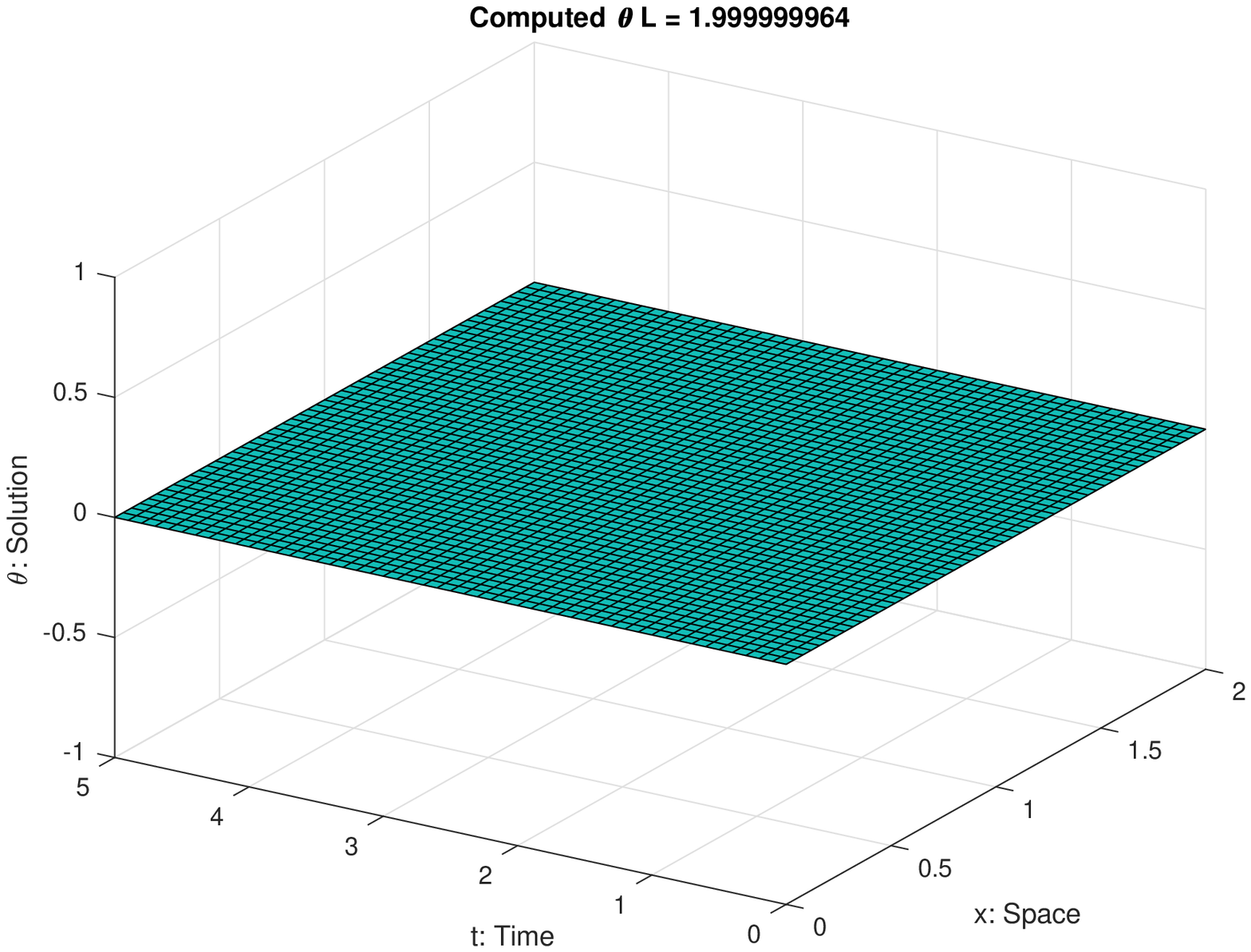}
\captionof{figure}{Burgers equation with heat effect with $(u_0,\theta_0)=(0,0)$ and $\eta\neq 0$ with one observation $u_x(0,t)$. The computed solution $\theta$.}
\label{BH_DNSolthC1}
\end{minipage}
\end{figure}

\begin{figure}[h!]
\begin{minipage}[t]{0.49\linewidth}
\noindent
\includegraphics[width=\linewidth]{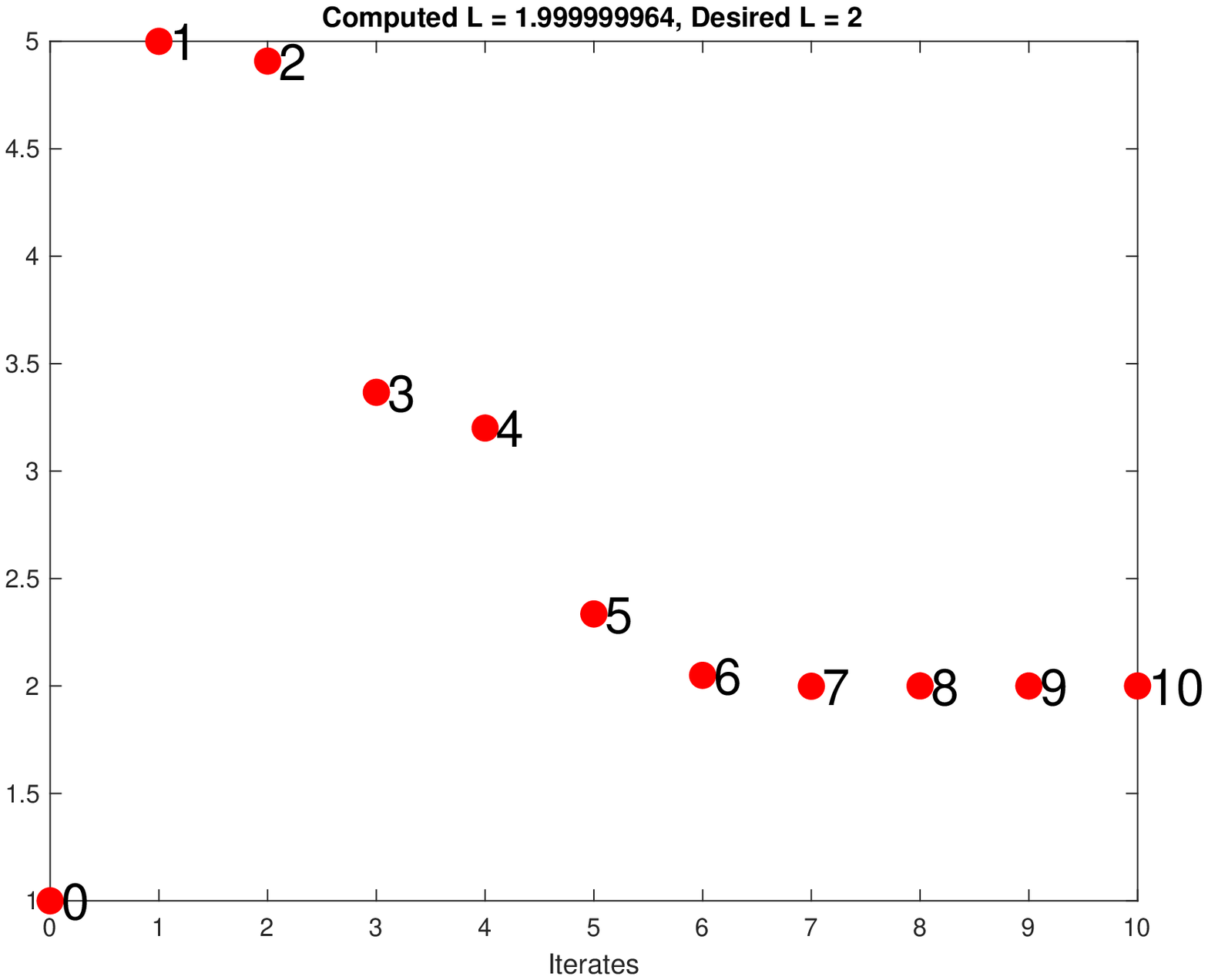}
\caption{Burgers equation with heat effect with $(u_0,\theta_0)=(0,0)$ and $\eta\neq 0$ with one observation $u_x(0,t)$. The iterates in \texttt{active-set} algorithm.}
\label{BH_DNiterC1}
\end{minipage}
\hfill
\begin{minipage}[t]{0.49\linewidth}
\noindent
\includegraphics[width=\linewidth]{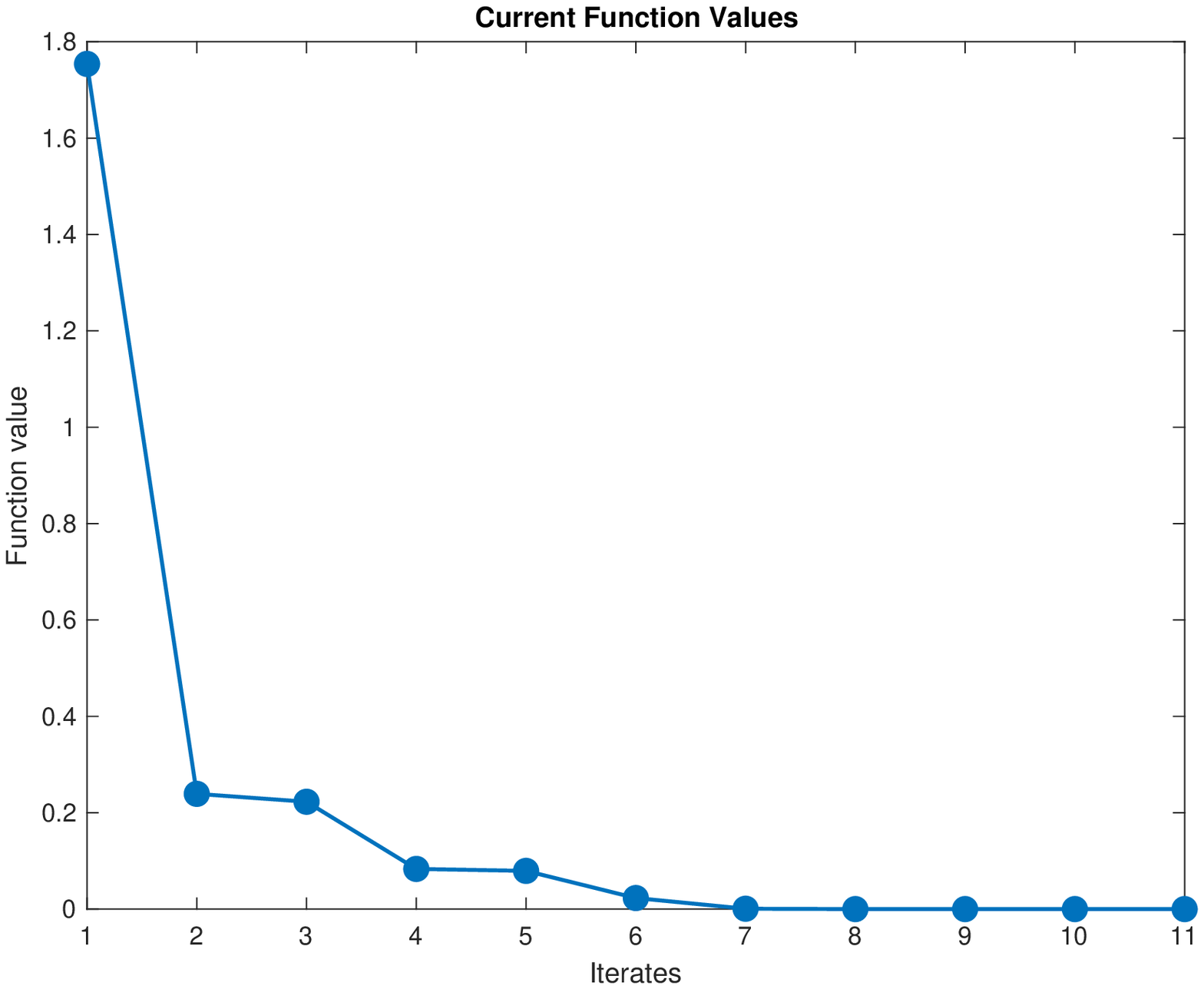}
\caption{Burgers equation with heat effect with $(u_0,\theta_0)=(0,0)$ and $\eta\neq 0$ with one observation $u_x(0,t)$. Evolution of the cost.}
\label{BH_DNfvC1}
\end{minipage}
\end{figure}

   The computed length is~$L_c=1.999999964$, the cost is $J(L_c)< 10^{-16}$ is reached in the iterate 10 of the optimization algorithm.
   The corresponding solution to~\eqref{pbBH} is displayed in~Figures~\ref{BH_DNSoluC1} and~\ref{BH_DNSolthC1}.  The evolution of the iterates and the cost in the minimization process in the absence of the random noise appear in~Figures~\ref{BH_DNiterC1} and~\ref{BH_DNfvC1}, respectively.
   
\

%

\noindent 
\textbf{Case 2.4: Burgers-heat system with $(u_0,\theta_0)\neq (0,0)$ and large $\eta$.}

 We take $T=5$,  $\eta(t) = 5\sin^3 t$ in~$(0,T)$ and~$u_0(x) = 0.1 x(2-x)$, $\theta_0(x) =0.1(1+x^2(x-3))$. Starting from $L_{i} = 1.4$, our goal is to recover the desired value of the length $L_d= 2$.

\begin{figure}[h!]
\begin{minipage}[t]{0.49\linewidth}
\vspace{0cm}
\includegraphics[width=\linewidth]{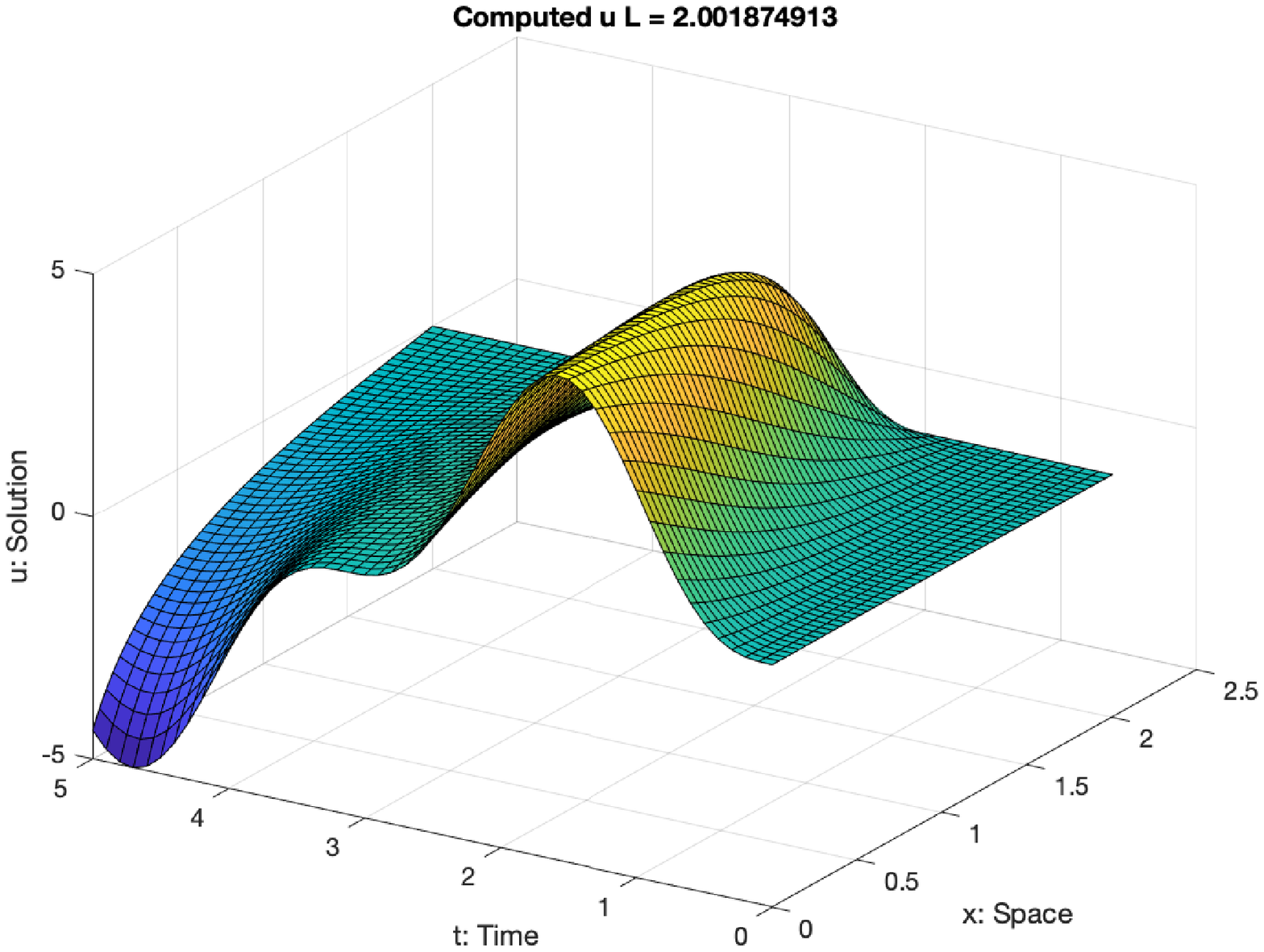}
\captionof{figure}{Burgers equation with heat effect with $(u_0,\theta_0)=(0,0)$ and $\eta\neq 0$ with one observation $u_x(0,t)$. The computed solution $u$.}
\label{BH_DNSoluC2}
\end{minipage}
\hfill
\begin{minipage}[t]{0.49\linewidth}
\vspace{0cm}
\includegraphics[width=\linewidth]{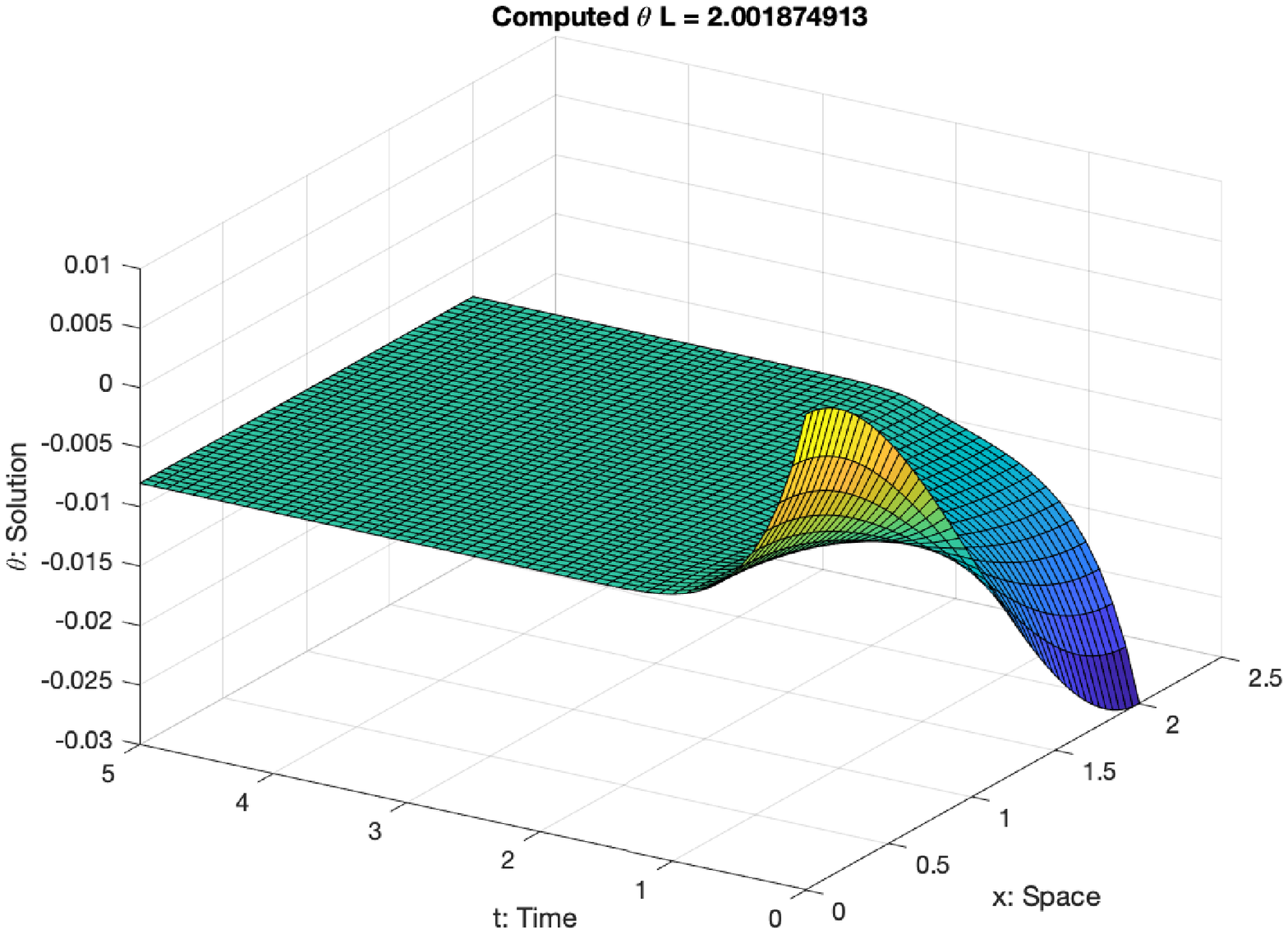}
\captionof{figure}{Burgers equation with heat effect with $(u_0,\theta_0)=(0,0)$ and $\eta\neq 0$ with one observation $u_x(0,t)$. The computed solution $\theta$.}
\label{BH_DNSolthC2}
\end{minipage}
\end{figure}

\begin{figure}[h!]
\begin{minipage}[t]{0.49\linewidth}
\noindent
\includegraphics[width=\linewidth]{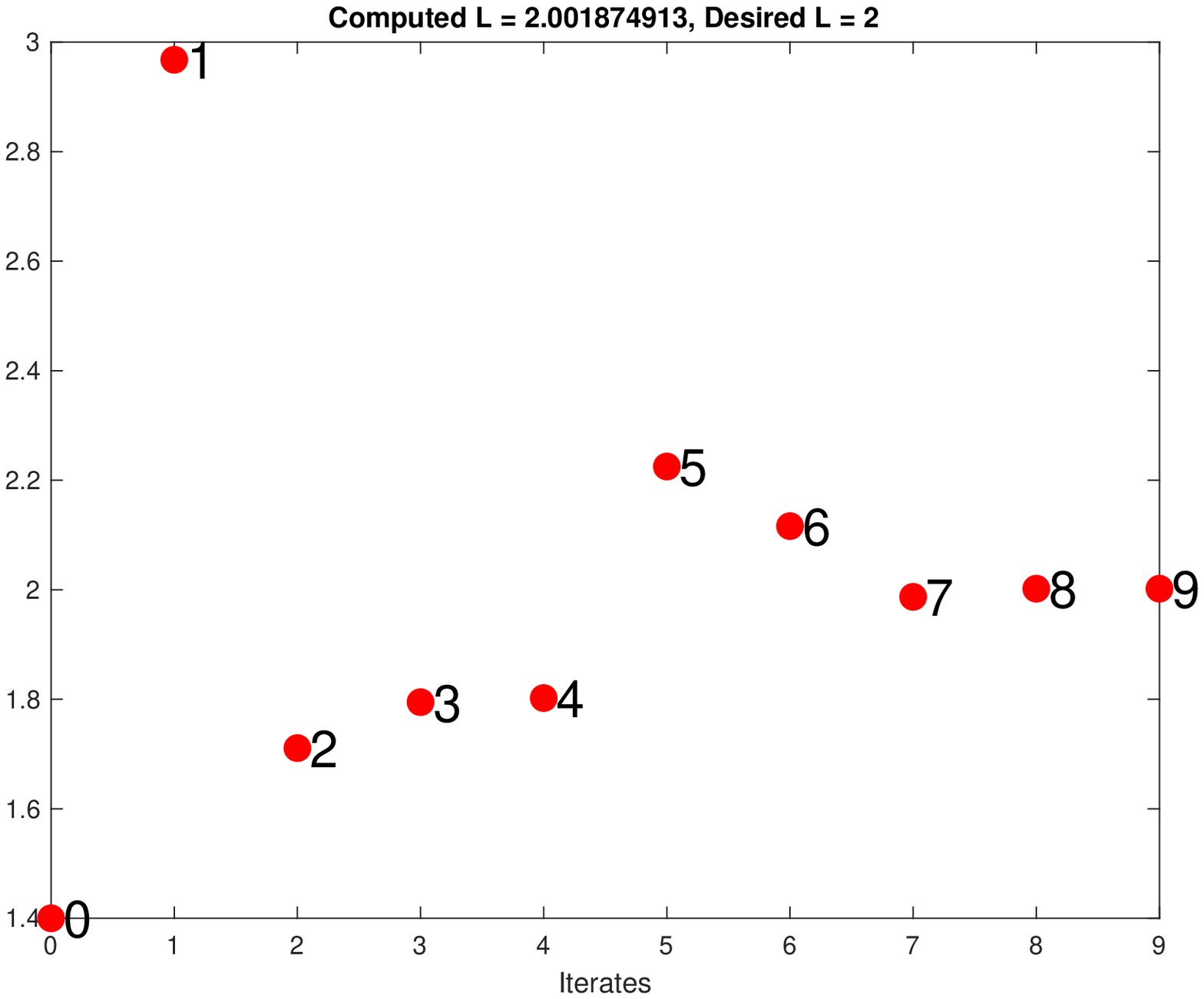}
\caption{Burgers equation with heat effect with $(u_0,\theta_0)=(0,0)$ and $\eta\neq 0$ with one observation $u_x(0,t)$. The iterates in \texttt{active-set} algorithm.}
\label{BH_DNiterC2}
\end{minipage}
\hfill
\begin{minipage}[t]{0.49\linewidth}
\noindent
\includegraphics[width=\linewidth]{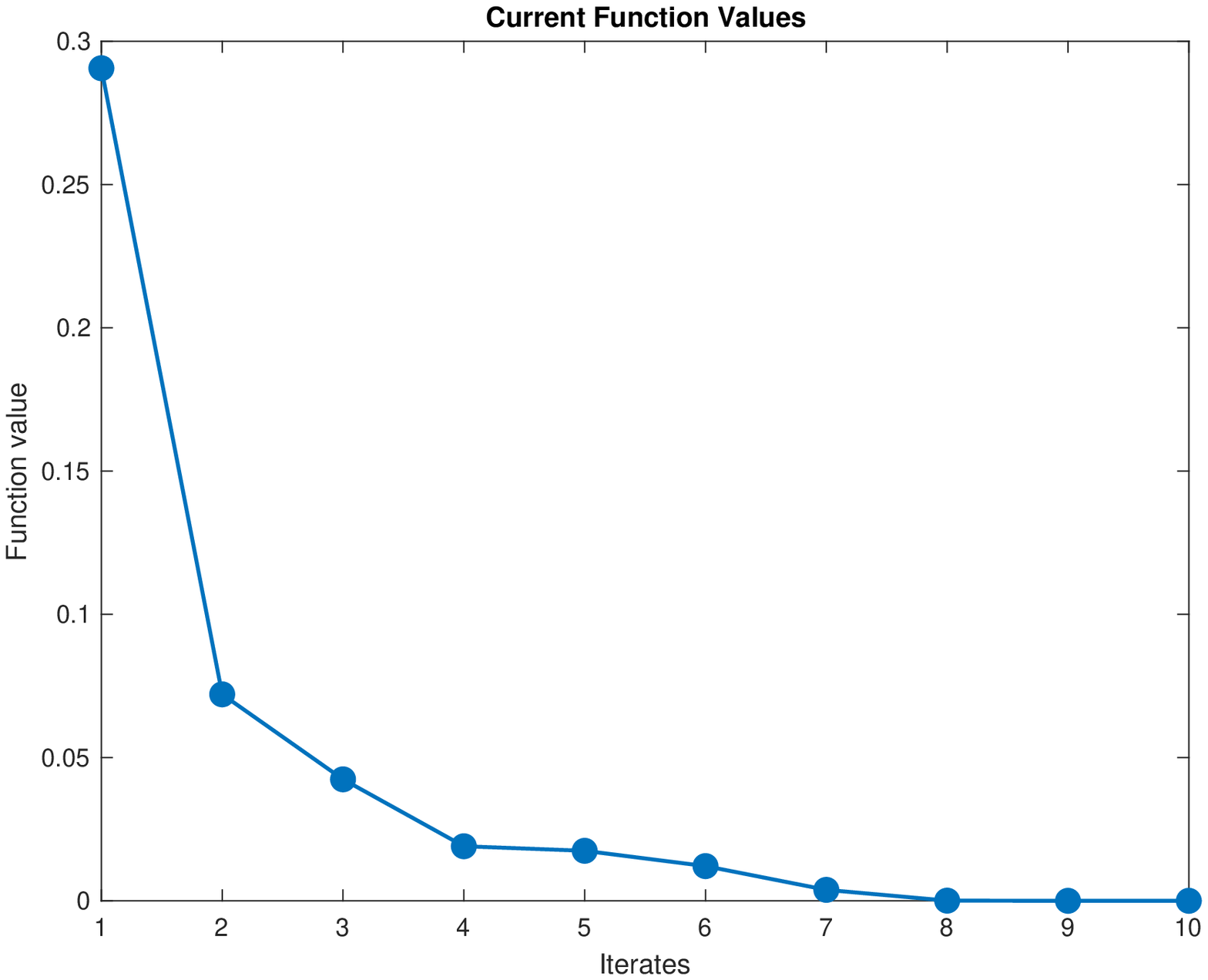}
\caption{Burgers equation with heat effect with $(u_0,\theta_0)=(0,0)$ and $\eta\neq 0$ with one observation $u_x(0,t)$. Evolution of the cost.}
\label{BH_DNfvC2}
\end{minipage}
\end{figure}

   The computed length is~$L_c=2.001874913$, the cost is $J(L_c)< 10^{-6}$ is reached in the iterate 9 of the optimization algorithm.
   The corresponding solution to~\eqref{pbBH} is displayed in~Figures~\ref{BH_DNSoluC2} and~\ref{BH_DNSolthC2}.  The evolution of the iterates and the cost in the minimization process in the absence of the random noise appear in~Figures~\ref{BH_DNiterC2} and~\ref{BH_DNfvC2}, respectively.
   

%
%

\section*{Acknowledgements}
The first author was supported by the Spanish Government's Ministry of Science, Innovation and Universities~(MICINN), under grant~PGC2018-094522-B-I00 and the Basque Government, under grant~IT12247-19.
The second and third authors were partially supported by~MICINN, under grant~MTM2016-76690-P. Lastly, the fourth author was supported by Grant-in-Aid for Scientific Research (S)
15H05740 of Japan Society for the Promotion of Science and by The National Natural Science Foundation of China
(no. 11771270, 91730303). This work was prepared with the support of the ``RUDN University Program 5-100".



\end{document}